%% file: main.tex
\documentclass[11pt,a4paper,twoside]{article}
\usepackage[T1]{fontenc}
\usepackage[utf8]{inputenc}
\usepackage{amssymb}
\input{macros.tex}
\begin{document}

\title{A Local  Discontinuous Galerkin Level Set Reinitialization with Subcell Stabilization on Unstructured Meshes}
\author[1]{A.~Karakus\footnote{Corresponding author. E-mail address: akarakus@metu.edu.tr} }
\author[2]{N.~Chalmers}
\author[3]{T.~Warburton}

\affil[1]{\textit{\small{Department of Mechanical Engineering, Middle East Technical University, Ankara, Turkey 06800}}}
\affil[2]{\textit{\small{Advanced Micro Devices Inc., Austin, TX}}}
\affil[3]{\textit{\small{Department of Mathematics, Virginia Tech, Blacksburg, VA 24061-0123}}}
\renewcommand\Authands{ and }
\date{\vspace{-5ex}}
\maketitle

\begin{abstract}
In this paper we consider a level set reinitialization technique based on a high-order, local discontinuous Galerkin  method on unstructured triangular meshes. A finite volume based subcell stabilization is used to improve the nonlinear stability of the method. Instead of the standard hyperbolic level set reinitialization, the flow of time Eikonal equation is discretized to construct an approximate signed distance function. Using the Eikonal equation removes the regularization parameter in the standard approach which allows more predictable behavior and faster convergence speeds around the interface. This makes our approach very efficient especially for banded level set formulations. A set of numerical experiments including both smooth and non-smooth interfaces indicate that the method experimentally achieves design order accuracy.  

\end{abstract}

\section{Introduction}
\input{introduction.tex}

\section{Formulation}
\input{formulation.tex}

\section{Results}
\input{results.tex}

\section{Conclusion}
\input{conclusion.tex}

\bibliographystyle{abbrv}
\bibliography{main}

\end{document}

%% file: macros.tex
\usepackage{lipsum}
\usepackage{authblk}
\usepackage[top=3cm, bottom=3cm, left=2cm, right=2cm]{geometry}
\usepackage{fancyhdr}
\usepackage{adjustbox}
 \usepackage[sort&compress, numbers, compress ]{natbib}
\usepackage[colorlinks,bookmarksopen,bookmarksnumbered,citecolor=red,urlcolor=red]{hyperref}
\usepackage{algorithm}
\usepackage{algorithmic}
\usepackage{boxedminipage}
\usepackage{comment}
\usepackage{eqparbox}
\usepackage{multirow}
\usepackage{caption}
\usepackage{subcaption}

\usepackage{amsmath}
\usepackage{amsfonts}
\usepackage{graphicx}

\usepackage{tikz}

\newlength{\commentindent}
\usepackage{array}
\newcolumntype{R}[1]{>{\raggedleft\let\newline\\\arraybackslash\hspace{0pt}}m{#1}}

\newcommand\BibTeX{{\rmfamily B\kern-.05em \textsc{i\kern-.025em b}\kern-.08em
T\kern-.1667em\lower.7ex\hbox{E}\kern-.125emX}}

\newcommand{\subc}[1]{\bar{#1}}
\newcommand{\EN}{\mathcal{E}}
\newcommand{\SEN}{\mathcal{S}}
\newcommand{\dE}{\partial \EN}

\newcommand{\grad}{\nabla}

\let\OLDthebibliography\thebibliography
\renewcommand\thebibliography[1]{
	\OLDthebibliography{#1}
	\setlength{\parskip}{0pt}
	\setlength{\itemsep}{0pt plus 0.3ex}
}

\renewenvironment{abstract}{%
	\begin{center}
		\begin{adjustbox}{minipage=0.90\textwidth,center} 
			\rule{\textwidth}{0.5pt}}
		{\par\noindent\rule{\textwidth}{0.5pt}  
		\end{adjustbox}
	\end{center} 
}

\makeatletter

\renewcommand\@maketitle{%
	\begin{adjustbox}{minipage=0.90\textwidth,center}
		\begin{center}
			\vskip 1em
			\let\footnote\thanks 
			{\LARGE \@title \par }
			\vskip 1.0em
			{\large \@author \par}
		\end{center}
	\end{adjustbox}
	\vskip 0.0em \par
}
\makeatother

%% file: introduction.tex
Level set (LS) methods \cite{osher_fronts_1988} are very popular to capture dynamic fronts in computational physics and engineering \cite{osher_level_2001,gibou_review_2018}. In a typical application, evolving the level set function in time often distorts the regularity of level set function. The reinitialization process replaces solution with signed distance function which satisfies the Eikonal equation $\lvert \nabla \phi \lvert = 1 $ by keeping the zero level set unchanged. 

In this study, we focus on the high-order, partial differential equation (PDE) based, reinitialization methods, and refer to Gibou et al. \cite{gibou_review_2018} for a recent review of a broad family of methods and their applications on various problems. The PDE-based level set reinitialization methods can be classified as pseudo-time first-order hyperbolic \cite{sussman_level_1994}, parabolic \cite{li_distance_2010},  and quasi-linear elliptic \cite{basting_minimization-based_2013} approaches. The former approach of Sussman \cite{sussman_level_1994} evolves the  Hamilton-Jacobi (HJ) equation in pseudotime to achieve signed distance function at steady state as given by
\begin{align}
\label{eq:RedistHyperbolic}
\frac{\partial \phi}{\partial t}+\text{sgn}_\alpha \left(\phi_{0}\right) \left(\lvert \nabla \phi \lvert-1\right) &=  0, \quad 
\phi\left(\mathbf{x},0\right) = \phi_{0}\left(\mathbf{x}\right),  
\end{align}
where $\text{sgn}_\alpha (\phi_{0})$ denotes $\alpha$ regularized sign term. The characteristics of this equation have unit speed and move from interface in the  normal direction without smoothing. However, the speed changes with  the $ \text{sgn}_{\alpha} (\phi_{0})$ term such that selection of the regularization  affects convergence rate, accuracy and stability of schemes \cite{sussman_improved_1998}. If $\alpha$ is selected as non-zero but too small a value, convergence speed of the solution increases with, however, sacrificed accuracy due to artificial movement of the interface. On the other hand, large values of $\alpha$ gives a profile as smooth as the initial data $\phi_0$ leading to slow convergence to steady state due to small characteristic speeds around the interface. Thus, in practice, $\alpha$ is usually chosen to be proportional to mesh size $h$, even though convergence property of the PDE is not explicitly known for general problems under different selections \cite{cheng_redistancing_2008}.

Most of the high-order level set reinitialization approaches for \eqref{eq:RedistHyperbolic}  are based on structured meshes utilizing  finite difference Essentially Non-Oscillatory (ENO) and Weighted ENO (WENO) schemes \cite{osher_high-order_1991,jiang_weighted_2000}. Although these schemes have been adapted to unstructured triangular grids \cite{zhang_high-order_2003,levy_central_2006} with highly increasing complexity, DG methods (we refer the reader to \cite{hesthaven_nodal_2008} and the references therein) have the advantages of easy implementation,  compactness, efficiency, and superior scalability. In \cite{grooss_level_2006}, the level set function is  reinitialized using a standard mixed DG method with an adaptive filtering  and stream line diffusion based stabilization, but both the filter strength and diffusion coefficients remained as parameters that may have to be tuned for a given problem. In \cite{karakus_gpu_2016,karakus_adaptive_2018}, the level set function was reinitialized on unstructured triangular/tetrahedral grids with a local DG \cite{yan_local_2011} space discretization paired with an artificial diffusion-type stabilization mechanism. To ensure accuracy in the stabilized elements, the mesh is refined to increase resolution and to decrease the diffusion added to system.  Zhang and Yue \cite{zhang_high-order_2019} recently introduced a gradient-based approach where the level set function is determined by a weighted local projection scheme adjusted to preserve the interface location for the elements involving the level set and used a HJ solver in conservative form for the remaining elements. 

The elliptic and parabolic reinitialization approaches are based on the minimization of an energy functional in strong form \cite{li_distance_2010,adams_parabolic_2019} or by directly applying penalty terms in an FEM framework to enforce boundary conditions at the interface \cite{basting_minimization-based_2013}. 
The elliptic reinitialization method was later analyzed using a DG space discretization for different potential functions and penalty parameters \cite{utz_interface-preserving_2017} motivated by the indicated issues of hyperbolic reinitialization using the DG framework \cite{mousavi_level_2014}. Applying Dirichlet boundary conditions on an immersed implicit interface requires special techniques such as Lagrange multipliers \cite{adams_high-order_2019} or constructing special quadrature rules for each element cut by the interface \cite{muller_highly_2013}. In fact, the smoothed sign term in the original level set formulation can be considered as a simplified implementation of Dirichlet conditions on the interface, and elliptic type reinitialization has similar problems to  interface preservation, i.e. the definition of a penalty parameter to enforce boundary conditions or the definition of high-order integration rules for the elements cut by interface involved in dynamic problems. Also, the computational cost of solving a quasi-linear elliptic problem can be prohibitively large for many practical applications. 

DG methods are well-suited for level set advection \cite{marchandise_quadrature-free_2006, karakus_gpu-accelerated_2016} or interface tracking in multiphase flows \cite{karakus_adaptive_2018, zhang_high-order_2019} due to their low numerical dissipation achieved by the use of local high order polynomial approximations. On the other hand, level set reinitialization is generally solved with a more robust, lower-order finite volume scheme \cite{fechter_discontinuous_2015} or is totally avoided using geometric distancing based on neighbor search algorithms \cite{marchandise_stabilized_2006,marchandise_stabilized_2007} even though the other parts of the solver benefit high-order DG discretizations.  This is probably because high-order DG methods for the HJ equation, like many other high-order methods, produce oscillations when the approximation space is inadequate to resolve the main features of the exact solution. Designing suitable stabilization mechanisms  for high-order discretizations of these problems is a challenging task on general problems and element types. Successful approaches for nonlinear stability are  limiting, which is based on reducing polynomial order near discontinuous regions \cite{cockburn_runge-kutta_1998}, high order WENO type polynomial reconstruction \cite{qiu_comparison_2005}, filtering high-frequency solution components, e.g. \cite{hesthaven_nodal_2008}, and artificial diffusion which relies on explicitly adding viscous terms to the governing equations \cite{persson_sub-cell_2006}. Recently, subcell-based limiters for DG \cite{huerta_simple_2012} have received attention since they act only on the smallest length scale within one cell to avoid the excessive numerical dissipation which might impact global solution. The \emph{a posteriori} subcell limiting process is often applied to tensor product elements since with the selection of Gauss-Lobatto interpolation nodes, together with a suitable sub-tessellation in the element, the DG and subcell FV flux computations match directly \cite{dumbser_posteriori_2014}.

In this work, we present a high-order local DG method for interface-preserving level set reinitialization on unstructured triangular meshes. Instead of using the standard hyperbolic reinitialization approach, we discretized flow  of time Eikonal equation which removes the dependency of the solution on the regularized sign term and improves convergence speed of solution. We also detail the efficiency of the method in the local set formulation where reinitialization is only needed in the vicinity of interface, as encountered in a multiphase flow simulation.  To stabilize the corresponding HJ equations, we designed \emph{a priori} subcell finite volume limiter on unstructured meshes that minimizes the additional operation count and connectivity information, and maximizes the amount of data gathered from each macro DG element. We  demonstrate that the proposed approach is robust, preserves $(N+1)$th order accuracy for smooth interfacial problems using $N$th order polynomial approximation space, and is suitable for fine-grain parallelism.   

The remainder of this paper is organized as follows. In Section 2, we present the mathematical formulation for the flow of the time Eikonal equation to compute signed distance functions, including the high-order local DG spatial discretizations and the construction of the level set function in time. Details of the subcell finite volume solver and \emph{a priori} error estimators are given in Section 3, which is followed by numerical validation test cases in Section 4. Finally, Section 5 is dedicated to concluding remarks and comments on future works. 

%% file: formulation.tex
The general approach for solving  the standard hyperbolic reinitialization \eqref{eq:RedistHyperbolic}   is to smear the signum term in a narrow band in the vicinity of the interface. Without regularization of the  sign function, characteristics of the equation are emanating from the interface in the  normal direction at unit speed. However, the speed changes with the $ \text{sgn}_{\alpha} (\phi_{0})$ term so that selection of the regularization affects the steady state solution and the convergence rate/accuracy of the scheme.  

Generally, regularized signum term is selected as 
\begin{equation}
     \text{sgn}_{\alpha} (\phi_{0})= \frac{\phi_{0}}{\sqrt{\phi_{0}^{2}+\alpha^2}}, 
\end{equation}
where $\alpha$ is the amplitude parameter chosen as a non-zero value related to the characteristic mesh size, i.e. $\alpha = O(h)$. Apparently, $\alpha$ should be large enough to ensure stability of the numerical discretization method and to prevent artificial movement of the interface. The stability issues become more prominent in high-order DG space discretizations \cite{mousavi_level_2014}. A careful calibration is required to get stable solutions \cite{karakus_gpu-accelerated_2016} under these conditions. In addition, a quick analysis of the characteristic equations for large values of $\alpha$ shows a source of inefficiency, since the characteristics that convey signed distance to space travel at speed $\text{sgn}_{\alpha} (\phi_{0})$, which is small close to the interface. This suggests that constructing the signed distance function close to the interface requires high amount of computation, and local level set reinitialization becomes more computationally demanding. 

Consider a uniformly continuous function, $\phi\left(\mathbf{x},t\right)$ which represents an interface, $\Gamma$ as
\begin{align*}
    \phi\left(\mathbf{x},0\right) &= 0, \quad \mathbf{x}\in \Gamma \\
    \phi\left(\mathbf{x},0\right) &> 0, \quad \mathbf{x}\in \Omega_p \\
    \phi\left(\mathbf{x},0\right) &< 0, \quad \mathbf{x}\in \Omega_m,
\end{align*}
where the domain $\Omega$ is partitioned as $\Omega_p \cup \Gamma \cup \Omega_m$, and $\mathbf{x}$ and $t$ are Cartesian coordinates and time, respectively. In \cite{osher_level_1993}, it is proved that  the interface evolves in time as $\phi\left(\mathbf{x},t\right) = 0$,  $\mathbf{x} \in \Gamma\left(t\right)$ for the given function $\phi\left(\mathbf{x},t\right)$ or $\frac{d\phi\left(\mathbf{x}\left(t\right),t\right)}{dt} = 0$. After a little manipulation, the relation can be written as
\begin{align} 
\label{eq:flow_of_time_Eikonal}
   \frac{\partial \phi}{\partial t} + \gamma\left(\mathbf{n}\right) \left\lvert \nabla \phi \right\rvert = 0,
\end{align}
where $\gamma\left(\mathbf{n}\right)$ is given as a function of the normal to level sets, $\mathbf{n}= \nabla\phi/\lvert \nabla \phi \rvert $. The velocity for  reinitialization can be constructed such that $\gamma\left(\mathbf{n}\right) = 1$, i.e. the flow is in the normal direction and has a speed of 1. Then \eqref{eq:flow_of_time_Eikonal} takes the following form, 
\begin{align}
\label{eq:Eikonal_single}
   \frac{\partial \phi}{\partial t} +\left\lvert \nabla \phi \right\rvert = 0.
\end{align}
Different from the original level set reinitialization approach, \eqref{eq:Eikonal_single} finds the first arrival times instead of distance level sets and is used to solve shape-from-shading  \cite{osher_level_1993} and redistancing \cite{cheng_redistancing_2008} problems successfully. Two flow fields are constructed as
\begin{equation}
\label{eq:redistance}
\begin{split}
    \frac{\partial u }{\partial t} + \left\lvert \nabla u \right\rvert &= 0, \quad u\left(\mathbf{x},0\right) = \phi_0\left(\mathbf{x}\right) \\
    \frac{\partial v }{\partial t} + \left\lvert \nabla v \right\rvert &= 0, \quad v\left(\mathbf{x},0\right) =- \phi_0\left(\mathbf{x}\right), 
\end{split}
\end{equation}
where $\phi_0$ is the initial value of the level set function. Therefore, for $\mathbf{x} \in \Omega_p$, the algorithm creates a flow such that the first arrival time of the front gives the signed distance in $\Omega_p$. The orientation of the front is reversed to find the arrival times for  $\mathbf{x} \in \Omega_m$.  Then the level set function is reconstructed accordingly as the arrival time, $T$ as follow
\begin{align}
\label{eq:time_recontruct}
    \phi\left(\mathbf{x}\right) = \left \{
    \begin{array}{l c c r}
        T   & \text{if } \phi_0\left( \mathbf{x}\right)>0  & \text{and} & u\left(\mathbf{x},T\right) = 0\\
        -T  & \text{if } \phi_0\left( \mathbf{x}\right)<0 & \text{and} & v\left(\mathbf{x},T\right) = 0\\
    \end{array} \right.
\end{align}
When the PDEs given in \eqref{eq:redistance} are discretized and the level set function is reconstructed in time using \eqref{eq:time_recontruct}, the accuracy of the numerical schemes used in time integration also plays an important role in the overall accuracy of the algorithm. We discuss the numerical discretizations we use along with the accuracy they provide in the following sections.

\subsection{Spatial Discretization}
We detail the spatial discretization of time-dependent Eikonal equation, \eqref{eq:redistance} in this section. We begin by partitioning the computational domain $\Omega$ into $K$ triangular elements $\EN^e$, $e=1,\ldots,K$, such that
\[
\Omega = \bigcup_{e=1}^K \EN^e.
\]
We denote the boundary of the element $\EN^e$ by $\dE^e$. We say that two elements, $\EN^{e+}$ and $\EN^{e-}$, are neighbours if they have a common face, that is $\dE^{e-} \cap \dE^{e+} \neq \emptyset$. We use $\mathbf{n} = \left(n_x, n_y\right)$ to denote the unit outward normal vector of $\dE$. 

We consider a finite element space on each element $\EN^e$, denoted $V_N^e  = \mathcal{P}_N(\EN^e)$ where $\mathcal{P}_N(\EN^e)$ is the space of polynomial functions of degree $N$ on element $\EN^e$. As a basis of the finite element spaces we take a set of $N_p = |V_n^e|$ Lagrange polynomials $\{l_n^e\}_{n=0}^{n=N_p}$, interpolating at the Warp $\&$ Blend nodes \cite{warburton_explicit_2006} mapped to the element $\EN^e$. Next, we define the polynomial approximation of the scalar level set field $\phi$ as 
\begin{align*}
    \phi^e &= \sum_{n=0}^{N_p} \phi^e_n l_n^e(\mathbf{x}), \text{ for all } \mathbf{x} = (x,y) \in \EN^e.
\end{align*}

The reinitialization equation, \eqref{eq:redistance}  can be written as a standard HJ equation in terms of generic field $\phi$ as $\phi = u$ or $\phi=v$. Considering  the first evolution equation $\phi=u$, since the second equation only uses flipped initial condition, we arrive to
\begin{equation}
\begin{aligned}
\label{eq.HJ_1}
\frac{\partial \phi}{\partial t}+H\left( \grad\phi\right)  =  0, \quad  \phi(\mathbf{x}, 0)  =\phi_{0},
\end{aligned}
\end{equation}
where the physical Hamiltonian, $H$ denotes $\lvert \grad \phi \rvert$. To solve the equation, the physical Hamiltonian is approximated by a monotone and consistent numerical Hamiltonian, $H\left( \grad \phi\right)\approx \bar{H}\left(\mathbf{p}, \mathbf{q} \right)$ such that $\mathbf{p}=\left[p_1,p_2\right]^T$ and $\mathbf{q}=\left[q_1,q_2\right]^T$ correspond to the vector of left and right upwind approximations of the gradient vector $\left(\phi_x, \phi_y\right)$, respectively:
\begin{equation*}
 \label{eq.directional_derivatives}
 \begin{split}
     p_1  - \phi_x^+ &= 0, \quad    p_2 - \phi_y^+ = 0 \\
    q_1  - \phi_x^- &= 0, \quad   q_2 - \phi_y^-  = 0.
 \end{split}
\end{equation*}
%
Then, we use the local Lax-Friedrichs type numerical Hamiltonian $\bar{H}\left(\mathbf{p}, \mathbf{q} \right)$ given as
\begin{equation}
\label{Eq.GodunovH}
\bar{H}\left(\mathbf{p}, \mathbf{q} \right)= H\left(\frac{p_1 + p_2}{2}, \frac{q_1 + q_2}{2}\right) -\alpha_1\left(p_1-p_2\right) -\alpha_2\left(q_1-q_2\right) .
\end{equation}
In the equation, $\alpha_1 = \max_{\mathbf{p}\in D}\left| \frac{\partial H\left(\mathbf{p}, \mathbf{q} \right)}{\partial \mathbf{p}}\right|$ and $\alpha_2 = \max_{\mathbf{q}\in D}\left| \frac{\partial H\left(\mathbf{p}, \mathbf{q} \right)}{\partial \mathbf{q}}\right|$ where domain $D$ is taken locally in the element $\EN^e$ as $D =\left[\min(p_1,p_2), \max(p_1,p_2)\right]$ or $D =\left[\min(q_1,q_2), \max(q_1,q_2)\right]$.  

Computing \eqref{Eq.GodunovH} requires accurate approximations of solution derivatives. We apply the  Local Discontinuous Galerkin (LDG) method \cite{yan_local_2011} due to its simplicity and minimal stencil size in unstructured triangular grids. 
Let $\mathbf{p}$ and $\mathbf{q}$ are left and right upwind approximation of gradient such that, 
\begin{align*}
    \mathbf{p} - \grad \phi = 0, \quad \mathbf{q} - \grad \phi = 0.
\end{align*}
Multiplying the equation by a test function $v \in V^e_N $, integrating over the element $\EN^e$, and performing integration by parts, we obtain the following weak variational form
\begin{equation}
\label{eq.WeakformEikonal}
\begin{split}
(v,\mathbf{p})_{\EN^e} = -(\nabla v,\phi)_{\EN^e} + (v,\phi_p^{*}\mathbf{n})_{\dE^{e}}\\
(v,\mathbf{q})_{\EN^e} = -(\nabla v,\phi)_{\EN^e} + (v,\phi_q^{*}\mathbf{n})_{\dE^{e}}.
\end{split}
\end{equation}
Here we have introduced the inner product $(u , v)_{\EN^e} $ to denote the integration of the product of $u$ and $v$ computed over the element $\EN^e$ and, analogously,  the inner product $(u , v)_{\dE^e}$ to denote the integration along the element boundary $\dE^e$. 

Due to the discontinuous approximation space, the flux functions $\phi$ are not uniquely defined in the boundary inner product and hence, it is replaced by a numerical flux function $\phi^*$ which depends on the local and neighboring trace values of $\phi$ along $\partial\EN^e$. On each element we denote the local trace values of $\phi^e$ as $\phi^{int}$ and the corresponding neighboring trace values as $\phi^{ext}$. Note that we will suppress the use of the $e$ superscript when it is clear which element is the local trace. Using this notation we choose as a numerical flux $\phi^*$ the alternating upwind fluxes as in \cite{yan_local_2011}, i.e. for $i=1,2$
\begin{equation}
\label{eq.upwind_flux}
\phi_{p}^{*}=
\begin{cases}
\begin{aligned}
\phi^{int} & \quad for \quad \mathbf{n}(i)\geqslant 0\\
\phi^{ext} & \quad for \quad  \mathbf{n}(i) < 0
\end{aligned}
\end{cases}    
\quad \text{and} \quad
\phi_{q}^{*}=
\begin{cases}
\begin{aligned}
\phi^{ext} & \quad for \quad \mathbf{n}(i) \geqslant 0\\
\phi^{int} & \quad for \quad\mathbf{n}(i) < 0
\end{aligned} 
\end{cases}  
\end{equation}

In order to write the discrete form of  \eqref{eq.WeakformEikonal} on the degrees of freedom of $\mathbf{p}=\left[p_1,p_2\right]^T$ and $\mathbf{q}=\left[q_1,q_2\right]^T$, we introduce the elemental mass $\mathcal{M}^e$, surface mass $\mathcal{M}^{ef}$, and stiffness operators $\mathcal{S}^e_x$ and $\mathcal{S}^e_x$  which are defined as follows
\begin{align}
    \mathcal{M}^e_{ij} = \left(l^e_i,l^e_j\right)_{\EN^e},  \label{eq:elMass}\quad
    \mathcal{M}^{ef}_{ij} = \left(l^e_i,l^e_j\right)_{\dE^{ef}}, 
\end{align}
\begin{align}
    (\mathcal{S}^e_x)_{ij} = \left(l_i^e,\frac{\partial l_j^e}{\partial x}\right)_{\EN^e},\label{eq:elStiff} \quad (\mathcal{S}^e_y)_{ij} = \left(l_i^e,\frac{\partial l_j^e}{\partial y}\right)_{\EN^e}.
\end{align}
Next, we define the weak elemental gradient operator $\boldsymbol{\mathcal{D}}^e = [\mathcal{D}^e_x, \mathcal{D}^e_y]^T$, as well as the lifting operators $\mathcal{L}^{ef}$, via
\begin{align}
    \mathcal{D}^e_x = (\mathcal{M}^e)^{-1}\mathcal{S}_x^T,\quad
    \mathcal{D}^e_y = (\mathcal{M}^e)^{-1}\mathcal{S}_y^T, \quad
    \mathcal{L}^{ef} = (\mathcal{M}^e)^{-1}\mathcal{M}^{ef}. \label{eq:elLift} 
\end{align}
Finally, for ease of notation we introduce the concatenation of the lift operators along each face, i.e., $\mathcal{L}^{e} = [\mathcal{L}^{e0}, \mathcal{L}^{e1}, \mathcal{L}^{e2}]$. Using these operators, fully discrete form of \eqref{eq.WeakformEikonal} can be written as

\begin{align*}
  \begin{split}
  p_1 = - \mathcal{D}_x^e \phi + \mathcal{L}^e \phi_{p,1}^*, \quad p_2 = - \mathcal{D}_y^e \phi + \mathcal{L}^e \phi_{p,2}^* \\
  q_1 = - \mathcal{D}_x^e \phi + \mathcal{L}^e \phi_{q,1}^*, \quad q_2 = - \mathcal{D}_y^e \phi + \mathcal{L}^e \phi_{q,2}^*,
\end{split}  
\end{align*}
which reduces performing excessive matrix-vector products in elemental computations.  In order to obtain a more unified expression between separate elements we introduce a mapping from each element $\EN^e$ to a reference element $\hat{\EN}$, on which we make use of reference operators. We take the reference element $\hat{\EN}$ to be the bi-unit triangle 
\[
\hat{\EN} = \left\{ -1\leq r,s,r+s\leq 1 \right\},
\]
and introduce the affine mapping $\Phi^e$ which maps $\EN^e$ to a reference triangle $\hat{\EN}$, i.e.,
\begin{equation*}
    \label{eq:operators1}
    \left(x,y\right) =\Phi^e\left(r,s\right), \quad \left(x,y\right)\in \EN^e,\ \left(r,s\right)\in\hat{\EN}
\end{equation*}
We denote the Jacobian of this mapping as
\begin{equation*}
    \label{eq:operators2}
    G^e = \begin{bmatrix}
    r_x & s_x \\
    r_y & s_y 
    \end{bmatrix},
\end{equation*}
and denote determinant of the Jacobian as $J^e = \det G^e$. We also define the surface scaling factor $J^{ef}$ which is defined as the determinant of the Jacobian $G^e$ restricted to the face $\partial \EN^{ef}$.  

Mapping each of the elemental operators defined in \eqref{eq:elMass}-\eqref{eq:elLift} to the reference element $\hat{\EN}$ we can write each of the elemental operators in terms of their reference versions and the geometric factors $G^e$, $J^e$, and $J^{ef}$ as follows
\begin{align*}
    \mathcal{M}^e = J^e \mathcal{M}, \quad 
    \boldsymbol{\mathcal{D}}^e = G^e \boldsymbol{\mathcal{D}},\quad \label{eq:elementOps}
    \mathcal{L}^{ef} = \frac{J^{ef}}{J^e}\mathcal{L}^f.
\end{align*}
Here $\mathcal{M}$, $\boldsymbol{\mathcal{D}} = [\mathcal{D}_r,\mathcal{D}_s]^T$, and $\mathcal{L}^f$ are the mass, derivative, and lifting operators defined on the reference element $\hat{\EN}$. Therefore, we can write the local DG scheme on each element using only these reference operators and the geometric data $G^e$, $J^e$, and $J^{ef}$. We also emphasize that volume terms of $p_1,q_1$ and $p_2,q_2$ share the same operations as seen in the discrete form. 

Finally, the semi-discrete form of the local DG scheme is written as,
\begin{equation}
\begin{aligned}
\label{eq.SemiDiscrete}
\frac{\partial \phi}{\partial t} = -\bar{H}\left( \mathbf{p}, \mathbf{q}\right).
\end{aligned}
\end{equation}
It is important to mention that when the solution is smooth $\mathbf{p}$ is very close to $\mathbf{q}$, while they differ significantly near the discontinuities. Thus, at discontinuous regions, ($\mathbf{p}$, $\mathbf{q}$) can capture the complete information of $\nabla \phi$. For a piecewise constant approximation, this scheme is monotone and converges to the entropy solution. However, stabilization is needed for higher order approximations. The use of subcell finite volume stabilization and  regularity detector is discussed in Section \ref{Subcell Stabilization}.

\subsection{Constructing Distance Function}
Overall accuracy in constructing the signed distance function depends on the time discretization scheme and the algorithm to find the first arrival time as given in \eqref{eq:time_recontruct}. For the time discretization, we use the low-storage fourth-order explicit Runge-Kutta scheme \cite{carpenter1994fourth} for the resulting ODEs. Then, the high-order determination of first arrival time of $u$ and $v$ remain to finish the overall discretization.

Any interpolation scheme can be used to compute an approximation of the time that $u$ and $v$ become zero. In order to achieve high order interpolation, ENO schemes offer desirable properties \cite{cheng_redistancing_2008}, especially if $u$ or $v$ has kinks in time. In this work, we consider a third order ENO interpolation that provides fourth order accurate approximation of the root when $u,v$ are smooth. We store the central time stencil for any grid point $\mathbf{x}_i$ as $\left\{u_i^{n-2}, u_i^{n-1},u_i^{n},u_i^{n+1},u_i^{n+2},u_i^{n+3}\right\}$ where $u_i^{n} = u\left(\mathbf{x}_i, t_n\right)$. Since roots of the ENO interpolating polynomial are always located in time interval $\left[t_n, t_{n+1}\right]$, we choose the initial stencil $E^{1} =\{t_i^{n},t_i^{n+1}\}$ which provides first order polynomial as 
\[p_1(t) = u[t_n,0] + u[t_n,1](t-t_n),\]
where any $u[t_j, k]$ is the $k$th order classical Newton divided difference for $k=0,1,2,3$ and $j=n-2,\ldots ,n+3$ given by
\[u\left[t_j,0\right] = u_i^{j}, \quad u\left[t_j,k+1\right] = \frac{u\left[j+1, k\right] - u\left[j, k\right]}{t_{j+k+1} - t_j}.\] 

Although selection of initial stencil $E^1$ is natural, since it is the closest to the root, the stencil $E^2$, which is required to construct a quadratic polynomial, is not fixed as either $E^2 = \{t_{n-1}, t_{n}, t_{n+1} \}$ or $
E^2 = \{t_{n}, t_{n+1}, t_{n+2} \}$ can be used to locate the root when $u$ is smooth in the interval. In this case, avoiding the candidate stencil which yields more a oscillatory interpolation is crucial to improve the accuracy. The ENO method \cite{harten_uniformly_1987} is to choose the stencil $S$ adaptively and automatically based on the local smoothness of the function $u$. The idea is to add the left or right neighbour point to the previous stencil each time step, depending on which one gives a smoother option using the highest degree term as indicator,
\begin{align*}
    E^{k+1} = \left\{\begin{array}{lcr}
        E^{k} \cup \{t_{l-1}\} & \text{ if } & \lvert u\left[t_{l-1}, k\right]\rvert \leq \left| u\left[t_l, k\right]\right|\\
        E^{k} \cup \{t_{l+k+1}\} & \text{ else} & 
        \end{array}\right.,
\end{align*}
where $l$ is the left most point in stencil $E^{k}$. Following this discussion, the $(k+1)$th order interpolating polynomial can be constructed as follows, 
\begin{align*}
    p_{k+1}(t) = p_{k}(t) + (t-E_1^{k})(t-E_2^{k})\ldots (t - E_{k+1}^{k})\left\{\begin{array}{lcc}
         u\left[t_{l-1}, k+1\right] & \text{ if } & \left| u\left[t_{l-1}, k+1\right]\right| \leq \left| u\left[t_l, k+1\right]\right|\\
         u\left[t_{l}, k+1\right]   & \text{ if } & \left| u\left[t_{l-1}, k+1\right]\right| > \left| u\left[t_l, k+1\right]\right|
    \end{array}\right..
\end{align*}
Here $E_l^{k}$ shows the $l$th entry of the stencil $E^{k}$. Then, classical root finding methods can be utilized to locate the roots of the polynomials. We use Newton’s method with initial guess of $(t^{n} + t^{n+1})/2$ , which gives a second order initial approximation with respect to time-step size. Our numerical tests demonstrated that only a couple Newton iterations allow us to find the root in targeted accuracy.

Close to the interface, there may not be enough values for an interpolation scheme to use. For example, if $n= 0$, then $u_i^0$ and $u_i^1$ are of different signs and we do not have sufficient number of history points. In order to use high order interpolation, we insert values from the approximation of $v$. Let $v_i^n$ denote the value computed by the scheme
to approximate $v(\mathbf{x}_i, t_n)$. Then we set $u_i^{-1} = -v_i^{1}$ and  $u_i^{-2} = -v_i^{2}$. This is because the level sets of $u$ and $v$ have initially the values of $\phi_0$ , and $u$ flows in its outward normal directions while $v$ flows in its inward normal directions.

\section{A Finite Volume Subcell Stabilization}
\label{Subcell Stabilization}
Effective stabilization of high-order DG discretizations of the reinitialization is critical to damp out high frequency solution components. Among the various techniques to prevent oscillations and to ensure stability in high-order discretizations, such as artificial diffusion and filtering, FV subcell limiting has favorable properties for nonlinear stability. Since all these mechanisms are designed to be integrated into DG schemes, detecting where the solution reaches its stability limits \emph{a priori} is crucial for effective solvers. In this work, we use a modal regularity detector introduced by Kl\"{o}ckner et.al. \cite{klockner_viscous_2011} as an improvement to the method proposed in \cite{persson_sub-cell_2006} by taking into account all the modes of expansion space for high-order DG solutions. This approach was adapted to 2D and 3D on triangular/tetrahedral elements in \cite{karakus_gpu_2016} where its effectiveness in terms of simplicity and marking of troubled elements in PDE based high-order level set reinitialization was shown.

The main idea of FV subcell stabilization is to replace the DG solution in troubled elements with a FV scheme on a subcell level, which is constructed by sub-tessellation of the macro DG element. Although FV subcell stabilization in a DG framework is well established in tensor product elements, as there is a natural balance between DG and FV data structures and communication patterns with suitable subdivisions \cite{sonntag_shock_2014}, the idea does not apply directly to simplices. Hence, special attention should be paid to designing efficient subcell limiting for triangular elements such that additional connectivity information storage size be minimized and information gathered from macro DG element should be maximized for FV solution.            

\begin{figure}[ht!]
	\begin{center}  
		\begin{subfigure}[]{0.3\textwidth}
			\includegraphics[width=\textwidth]{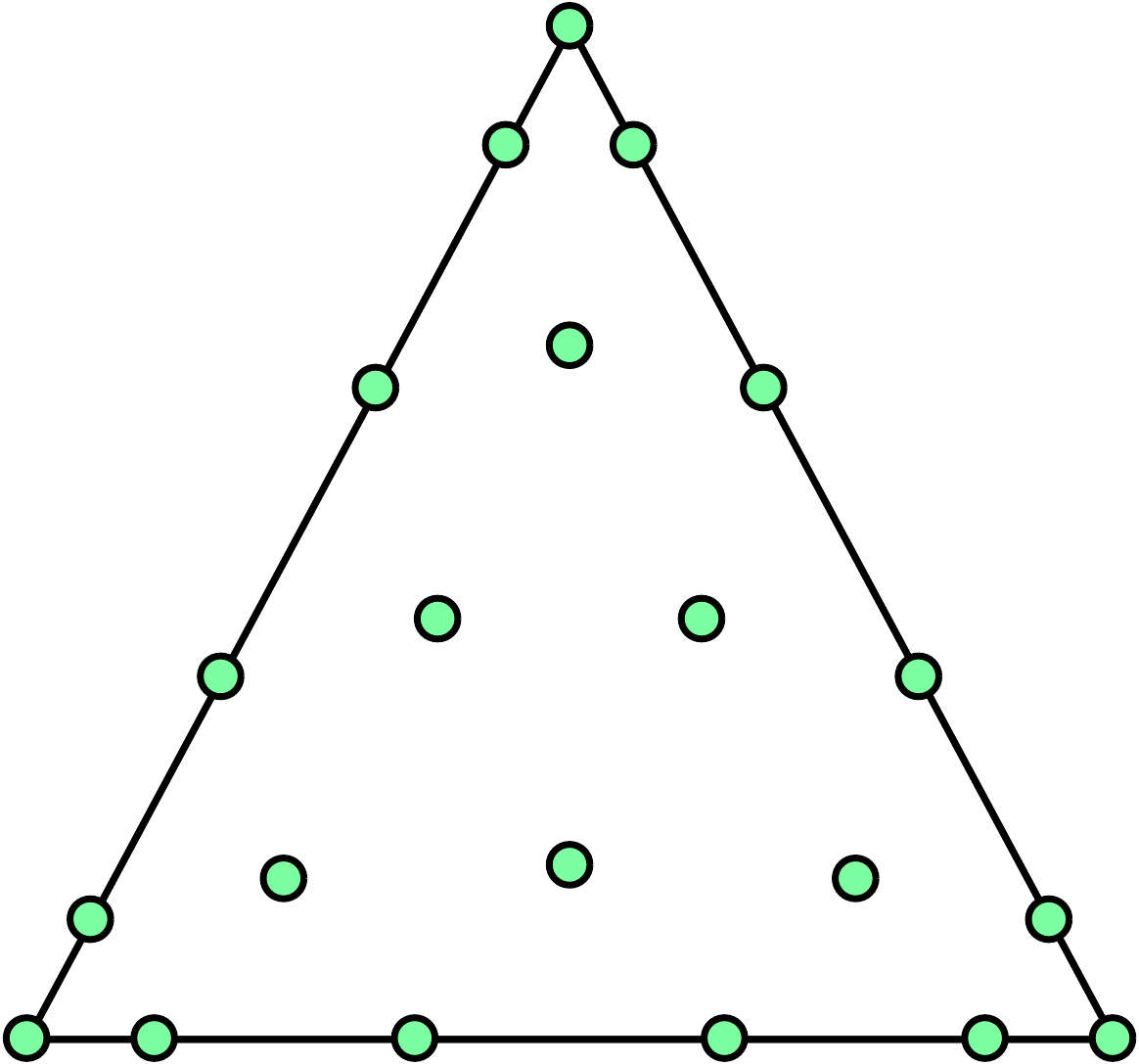}
			\label{Fig:WB_nodes_N5}
			\caption{}
		\end{subfigure} 
		~
		\begin{subfigure}[]{0.3\textwidth}
			\includegraphics[width=\textwidth]{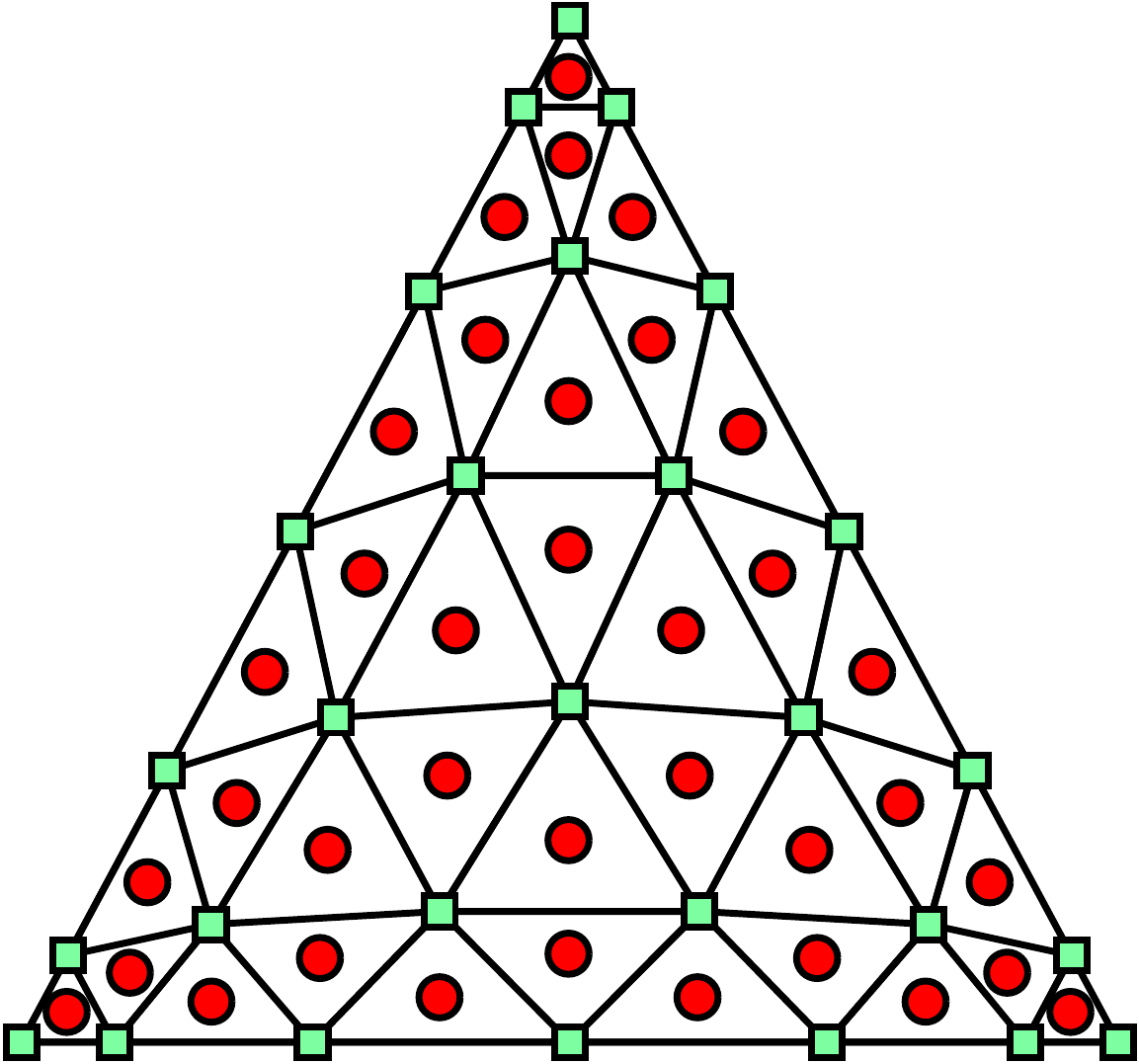} 
			\label{Fig:Subcell_WB_N6}
			\caption{}        
		\end{subfigure} 
	\end{center}
	\caption{Interpolation nodes and subcells used for \emph{a priori} finite volume limiting.  Warp$\&$Blend interpolation nodes \cite{warburton_explicit_2006} for $N=5$ (a) and subcells and their centroids created on the same node distribution for $N_s = (N+1)^2$ (b). }
	\label{Fig:WBNodes_and_Subcells}
\end{figure}

The space of $N$th order Lagrange polynomial functions on a triangular element can be constructed using  $N_p = |V_n^e|$ nodes, which has $N+1$ nodes on each face. The same polynomial space can be obtained from cell averages of  at least  $N_s = (N+1)^2$ subcells. In other words, any function represented by DG finite element space can be represented equivalently by a collection of piecewise constant finite volume cell data on a subgrid for $N_s\geq (N+1)^2$.   Although the number of FV subcells and the type of minor tessellation can be obtained in different ways, we choose  $N_s=(N+1)^2$ which is constructed at  $(N+1)$th order  Warp $\&$ Blend nodes \cite{warburton_explicit_2006} as shown in the figure \ref{Fig:WBNodes_and_Subcells} for $N=5$. The motivation of our particular choice is both to keep computational complexity of DG and FV cell workloads as close as possible and to preserve the subscale resolution of $N$th order DG method the same with FV characteristic length scale of $h/(N+1)$ \cite{ainsworth_dispersive_2004}.  

Let $\subc{\phi}\left(\mathbf{x},t\right)$ represent $\phi\left(\mathbf{x},t\right)$ with the same nominal accuracy, i.e. $N_s\geq (N+1)^2$,  in terms of a set of piecewise constant subcell averages $\subc{\phi}_i$ for $i=1\ldots N_s$. Then the $\subc{\phi}_i^e$ are computed as the $L_2$ projection of $\phi^e$ onto the space of piecewise constant functions on subcell $\mathcal{S}_i^e$ on macro element $\EN^e$ given by,
\begin{equation}
    \subc{\phi}_i^e = \frac{1}{\lvert \mathcal{S}_i^e\rvert} (\phi^e, \mathbf{1})_{\SEN_i^e} =  \frac{1}{\lvert \mathcal{S}_i^e\rvert} (l_n^e, \mathbf{1})_{\SEN_i^e} \phi_n^e,  \quad \forall \mathcal{S}_i^e \in \EN^e, n=1\ldots N_p.
\end{equation}
Since the integrand only includes the Lagrange interpolating polynomials over the subcells, it can be written on the macro reference element and its subdivision as follow,
\begin{equation}
   \mathcal{P}_{i,j} =  \frac{1}{\lvert\hat{\mathcal{S}}_i\rvert} \left(l_n,\mathbf{1}\right)_{\hat{\mathcal{S}}_i} \quad \forall \hat{\mathcal{S}}_i \in \hat{\EN}, n=1\ldots N_p.
\end{equation}
This defines the projection operator $\mathcal{P}$ that gives elemental mean values at subcells via $ \subc{\phi}^e = \mathcal{P} \phi^e$. We also define the projection operator restricted to the face, $\partial \hat{\EN}^{f}$,
\begin{equation}
   \mathcal{P}^f_{i,j} =  \frac{1}{\lvert\hat{\mathcal{S}}_i^f\rvert} \left(l_n,\mathbf{1}\right)_{\hat{\mathcal{S}^f}_i} \quad \forall \hat{\mathcal{S}^f}_i \in \hat{\EN}^f, \quad i=1\ldots (N+2),  j=1\ldots (N+1).
\end{equation}
$\mathcal{P}^f$ plays an important role in conservative couplings of a DG element with a FV subcell as piecewise constant mean values at faces can be obtained by $\bar{\phi}^{e,f} = \mathcal{P}^f \phi^{e,f}$.   

Reconstruction of a high-order polynomial representation from subcell averages  obviously leads an over-determined system as $N_s>N_p$ by construction. Hence this recovery problem can be solved in constrained least squares sense, where the constraint preserves the mean value of the high-order polynomial reconstruction on macro element,$ \left(\bar{\phi}, \mathbf{1}\right)_{\EN^e} = \left(\phi, \mathbf{1}\right)_{\EN^e}$. Following this relation, the reconstruction operator, $\mathcal{R}$ becomes pseudo-inverse of projection operator, $\mathcal{R}$ such that $\mathcal{R} \mathcal{P} = \mathcal{I}$, where $\mathcal{I}$ denotes the identity matrix. Similarly, reconstruction from face mean values satisfies $\mathcal{R}^f \mathcal{P}^f = \mathcal{I}$.  

After obtaining cell averages on the subgrid, any stable finite volume scheme can be used to evolve \eqref{eq:redistance}. Since we choose the minimum required number of subcells for minor tessellation to get an efficient scheme, the quality of the FV update is critical to achieve accurate solutions in troubled cells. The accuracy of a finite volume scheme is highly dependent on the approximation quality of gradients. Obviously, the simplest choice is to use a first order scheme which assumes a constant field or zero gradient over the subcells, but this results in a highly dissipative scheme even in subcell resolution on unstructured grids.  To obtain second-order accuracy in the FV subcells, a gradient scheme is needed to be at least first-order in general unstructured grids. Although the Green-Gauss gradient reconstruction is a computationally attractive technique, it requires special conditions on unstructured grids to achieve first-order accuracy in gradient. Alternatively, a least-squares gradient reconstruction, which minimizes the magnitude of gradient, provides first-order accuracy unconditionally on all grids but requires that the obtained linear system is not singular. 

\begin{figure}
\begin{center}
    \begin{subfigure}[b]{0.35\textwidth}
       \includegraphics[width=\textwidth]{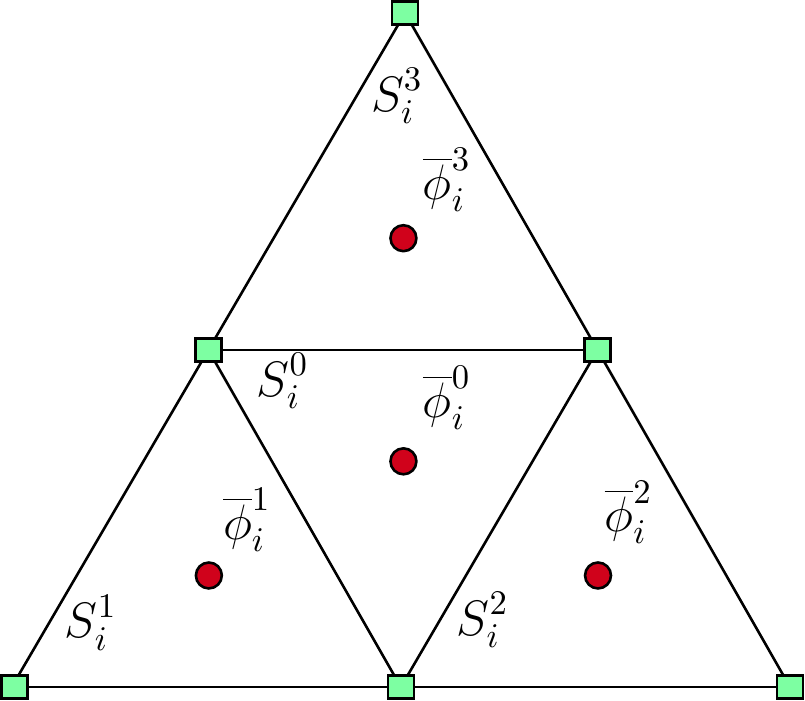} 
       \caption{}
    \end{subfigure}
    \begin{subfigure}[b]{0.4\textwidth}
       \includegraphics[width=\textwidth]{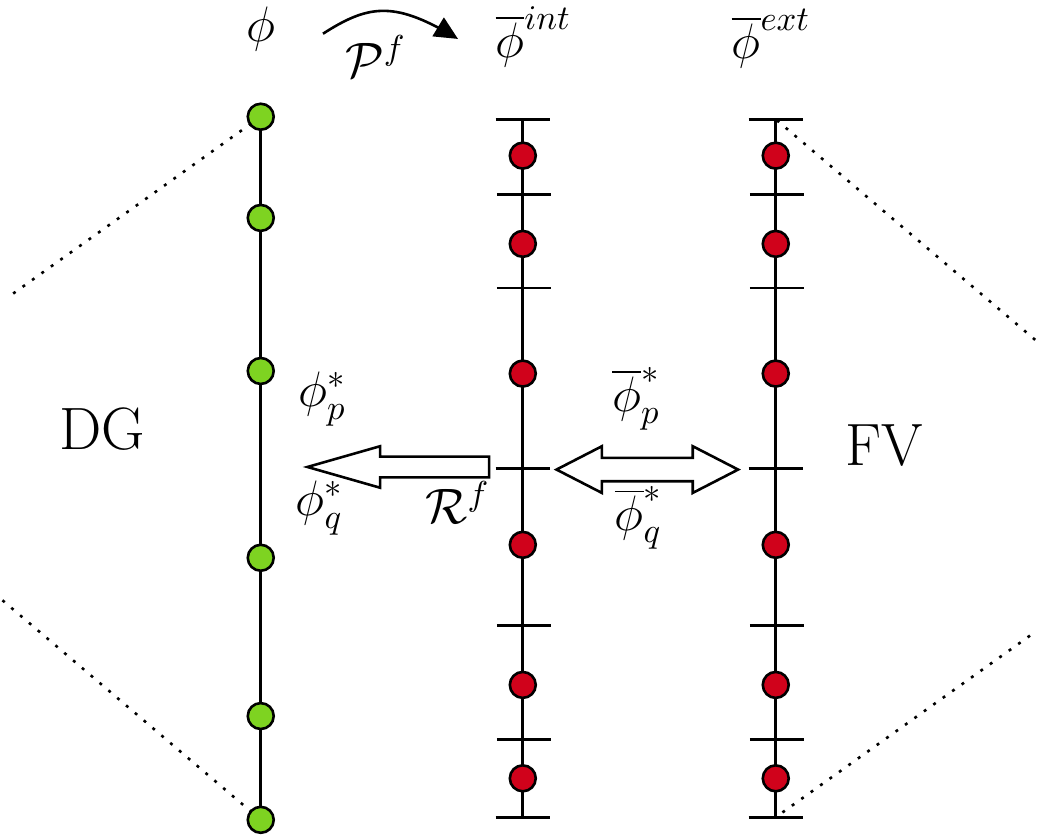} 
       \caption{}
    \end{subfigure}
\end{center}
\caption{(a) Finite volume stencil and data structure for WENO reconstruction and (b) computation of the flux function on a face of a DG element which is connected to FV elements.}
\label{fig:WENO_stencil_coupling}
\end{figure}

In this study, we utilize a second-order WENO approach \cite{friedrich_weighted_1998} to reconstruct gradients on subcell centers which results in a first-order accurate gradient approximation on general grids and avoids singularity problems with proper selection of nonlinear weights. 
For second-order accuracy, the degree of the required polynomial approximation of the gradient is $1$, which  can be constructed from mean values with stencil size of $3$. The polynomials $p_j$ for $j=1,2,3$ are constructed for the subcell $S_i^0$ using the compact stencils $\{S_i^0, S_i^1, S_i^2\}, \{S_i^0, S_i^2, S_i^3\}$ and  $\{S_i^0, S_i^3, S_i^1\}$ as illustrated in the figure \ref{fig:WENO_stencil_coupling}. We use the standard notation for the reconstruction polynomials as $p_j = a_j^0 + a_j^1 \left(x - \bar{x}_i^0\right) + a_j^2\left(y - \bar{y}_i^0\right) $ where $(\bar{x}, \bar{y})$ represents center coordinates of the indicated subcell. The unknowns $a$ are uniquely solvable for mean values $\bar{\phi}$ for the the given admissible stencils. By construction, $a_j^0 = \bar{\phi}_i^0$ for all polynomials and  $[a_j^1,a_j^2]$ are the components of gradient at $S_i^0$ computed on the  stencil $j$. Then, the resulting linear system for the first stencil can be written as 
\begin{equation*}
    \begin{bmatrix}
    \bar{x}_i^1-\bar{x}_i^0 & \bar{y}_i^1-\bar{y}_i^0\\
    \bar{x}_i^2-\bar{x}_i^0 & \bar{y}_i^2-\bar{y}_i^0
    \end{bmatrix}
    \begin{bmatrix}
    a_1^1\\
    a_1^2
    \end{bmatrix} = \begin{bmatrix}
    \bar{\phi}_i^1 - \phi_i^0\\
    \bar{\phi}_i^2 - \phi_i^0\\
    \end{bmatrix},
\end{equation*}
and similar for the stencils $2$ and $3$. For the subcells connected on DG elements, required stencils are constructed considering the corresponding face mean values of connected DG elements instead of  using the mean values of subcells created on dummy minor tessellation of neighbor DG elements.      

In WENO schemes, the reconstruction polynomial $p\left(\mathbf{x}\right)$ is obtained as a weighted sum of all polynomials, $p\left(\mathbf{x}\right)=\sum_{j} \omega_j p_j\left(\mathbf{x}\right)$.  The weights $w_j$ are chosen depending on smoothness of $p_j$ to minimize oscillations.  We use a normalized $L_2$ norm on the first derivative of $p_j\left(\mathbf{x}\right)$ to detect how much it oscillates 
\begin{equation*}
    \gamma_j = \left| \mathcal{S}_i \right|^{-2} \left( \lvert \grad p_j\left( \mathbf{x}\right)\rvert ^2,1 \right)_{\mathcal{S}_i}^\frac{1}{2}
\end{equation*}
which reduces to the simpler form of  $\gamma_j = \left| \mathcal{S}_i \right|^{-1} \left( (a_j^1)^2 + (a_j^2)^2 \right)^{\frac{1}{2}}$ for second-order reconstruction polynomials. The weights are then computed as 
\begin{equation*}
    \omega_i = \frac{(\epsilon + \gamma_i)^{-r}}{\sum_{j=1}^3 \left( \epsilon + \gamma_j \right)^{-r}}
\end{equation*}
where $r$ is a positive integer selected as $4$ and $\epsilon$ is a small positive number that we take as $10^{-6}$. Our numerical tests show that results are not overly sensitive to selection of $r$ and $\epsilon$. 

$p\left(\mathbf{x}\right)$ is then used to obtain face values on the subcells, $\bar{\phi}_i^f$ in $S_i$ which follows the discretization of the flow of the time Eikonal equation \eqref{eq.HJ_1}. With the piecewise constant approximation space, volume terms in the LDG discretization described in \eqref{eq.WeakformEikonal} vanishes and the scheme degenerates to 
\begin{equation}
  \label{eq.WeakformEikonalFV}
\begin{split}
\bar{\mathbf{p}}_i = \frac{1}{\lvert \SEN_i \rvert} \sum_{f=1}^{N_f} \bar{\phi}_p^{*}\lvert\SEN_i^f\rvert \mathbf{n}^f, \quad
\bar{\mathbf{q}}_i = \frac{1}{\lvert \SEN_i \rvert} \sum_{f=1}^{N_f} \bar{\phi}_q^{*}\lvert\SEN_i^f\rvert \mathbf{n}^f,
\end{split}  
\end{equation}
where $\bar{\mathbf{p}}_i = [\bar{p}_1, \bar{p}_2]^T$ and $\bar{\mathbf{q}}_i = [\bar{q}_1, \bar{q}_2]^T$ correspond to the vector of left and right upwind approximations of the gradient vector $[\bar{\phi}_x, \bar{\phi}_y]^T$ on $\mathcal{S}_i$. The alternating upwind fluxes,$\bar{\phi}_p^*$ and $\bar{\phi}_q^*$ on subcell faces $\mathcal{S}_i^f$ are defined similarly to the \eqref{eq.upwind_flux}. Finally, with the computation of derivatives for internal and boundary elements, we arrive to the semi-discrete form by using a LLF numerical Hamiltonian,
\begin{equation}
\begin{aligned}
\label{eq.SemiDiscreteFV}
\frac{\partial \bar{\phi}}{\partial t} = -\bar{H}\left( \bar{\mathbf{p}}, \bar{\mathbf{q}}\right).
\end{aligned}
\end{equation}
which completes the spatial discretization on subcells.

The degenerated scheme simply leads to a block-structured upwind FV scheme which offers inherent coupling of DG and FV elements through the same numerical fluxes. The numerical flux for a face of DG element is computed  on FV side using face projection operator, $\mathcal{P}^f$ as shown in Figure \ref{fig:WENO_stencil_coupling} (b). After this operation, external and internal trace values coincide at the FV element's  degrees of freedom where upwind fluxes are evaluated. These fluxes are used to update the FV subcells directly and are reconstructed back to DG degrees of freedom to evolve the solutions in the DG elements i.e. $\phi_p^* = \mathcal{R}^f\bar{\phi}_p^*$ and $\phi_q^* = \mathcal{R}^f\bar{\phi}_q^*$. 

%% file: results.tex
\label{Results}
In this section, we show the accuracy in local and global LS formulations, interface preservation, and long-term stability 
of the presented algorithm. Reinitialization tests are solved for both smooth and non-smooth interface problems on 2D unstructured triangular grids. To evaluate the accuracy of the  numerical scheme, we define the following $  L_{\infty} $ and $ L_{2} $ norms of error according to the following relations,
\begin{equation}
\label{eq.ErrorNorms1}
\begin{aligned}
L_{\infty}  &= \max_{\EN^e\in \Gamma_\epsilon} \lvert\phi^e-\phi^{exact}\lvert,\quad
L_{2}       &= \sum(\phi^e-\phi^{exact}, \phi^e-\phi^{exact})_{\EN^e\in \Gamma_\epsilon}
\end{aligned}
\end{equation}
where $\epsilon $ is the predefined band thickness and $\Gamma_\epsilon$ is the $\epsilon $  neighbourhood of the interface or whole computational domain depending on the location where the error is computed.    $\phi^{exact}$ denotes the exact solution or a very accurate approximation of the exact solution if it is not explicitly known. To assess the interface preservation of the method, we introduced the averaged $L_1$ error which measures displacement of the interface
\begin{equation}
\label{eq.ErrorNorms2}
\begin{aligned}
L_{1}       = \frac{1}{L_{\Gamma}}(\mathcal{H}_h \left( \phi^e\right)-\mathcal{H}_h \left( \phi^{exact}\right),\mathbf{1})_{\EN^e},
\end{aligned}
\end{equation}
where $L_{\Gamma}$ is length of the interface and $\mathcal{H}_h$ is a smoothed Heaviside function defined as $\mathcal{H}_h \left(\phi\right)=0.5\left(1+ \tanh(\frac{\pi\phi}{h})\right)$, $h$ being the characteristic length of mesh. 

In all the test cases, unless explicitly stated otherwise, we use $N_s=N+1$ for the subcells and the low storage, fourth-order explicit Runge-Kutta \cite{carpenter1994fourth} (LSERK) time discretization with unit CFL number. 

\subsection{Smooth Interfaces}
\label{Sec.Smooth Interfaces}
In the first test case, the reinitialization problem proposed by \cite{sussman_improved_1998} is solved to show convergence rate and interface preservation of the algorithm. The interface of interest is a circle centered at the origin with radius of $ 1.0 $. The signed distance function in a computational domain of $ [-2,2]^2 $  is perturbed into the following initial level set function
\begin{equation}
\label{Eq.Smooth_Circle_IC}
\phi_{0}=((x-1)^{2}+(y-1)^{2}+0.1) (\sqrt{x^{2}+y^{2}}-1),
\end{equation}
to create a smooth function with widely changing gradients in the vicinity of the interface. Computations are performed on the successively refined meshes by uniformly dividing the initial coarse level grid which is constructed by triangles having edge length of  $h=0.4$ on the boundary of domain, $\Omega$ and element number of $K=240$.  
\begin{figure}[ht!]
	\begin{center}
		\begin{subfigure}[b]{0.3\textwidth}
				\includegraphics[width=\textwidth]{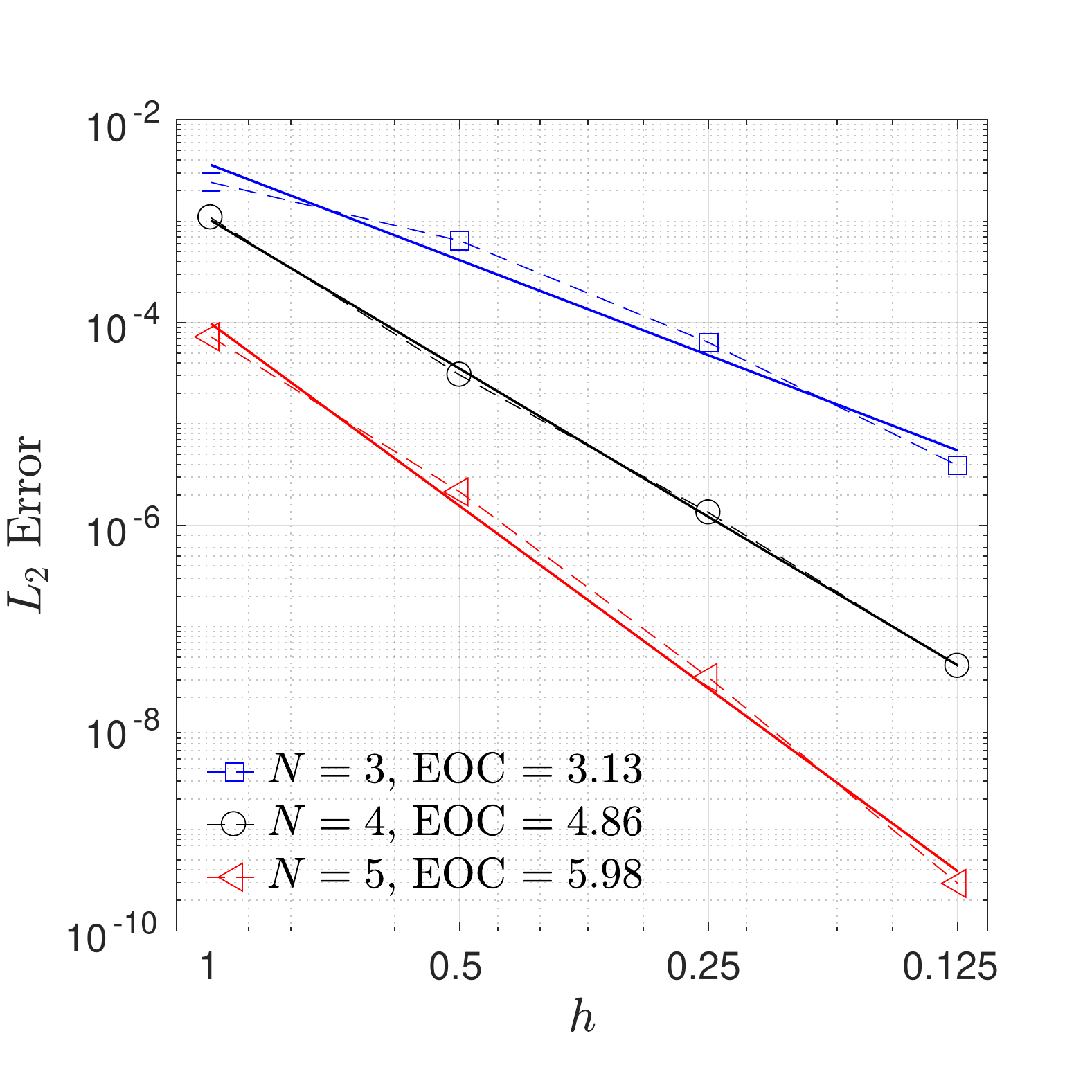}  
		\end{subfigure}
		        ~
		        \begin{subfigure}[b]{0.3\textwidth}
		           \includegraphics[width=\textwidth]{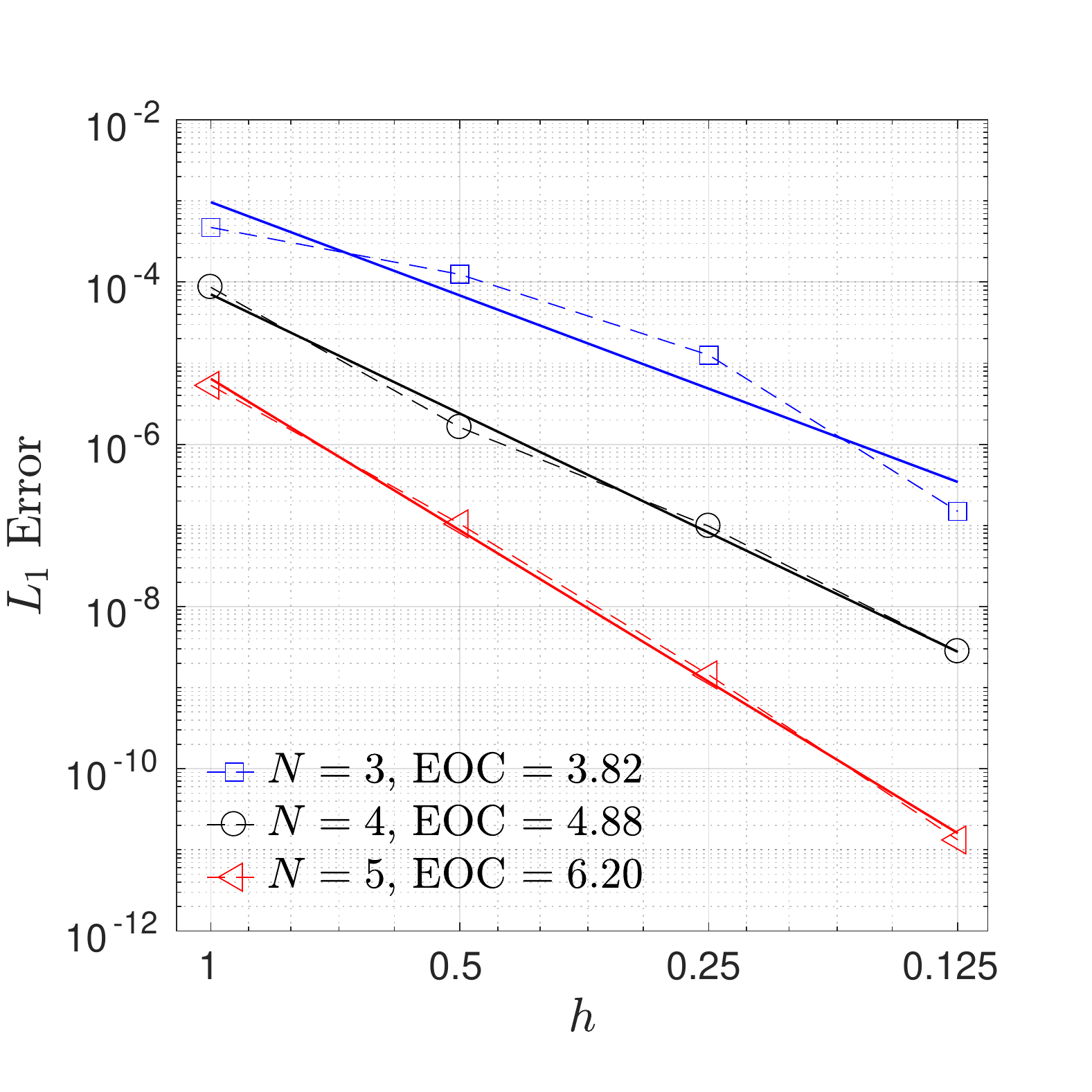}
		       \end{subfigure}
		        ~
		        \begin{subfigure}[b]{0.3\textwidth}
		           \includegraphics[width=\textwidth]{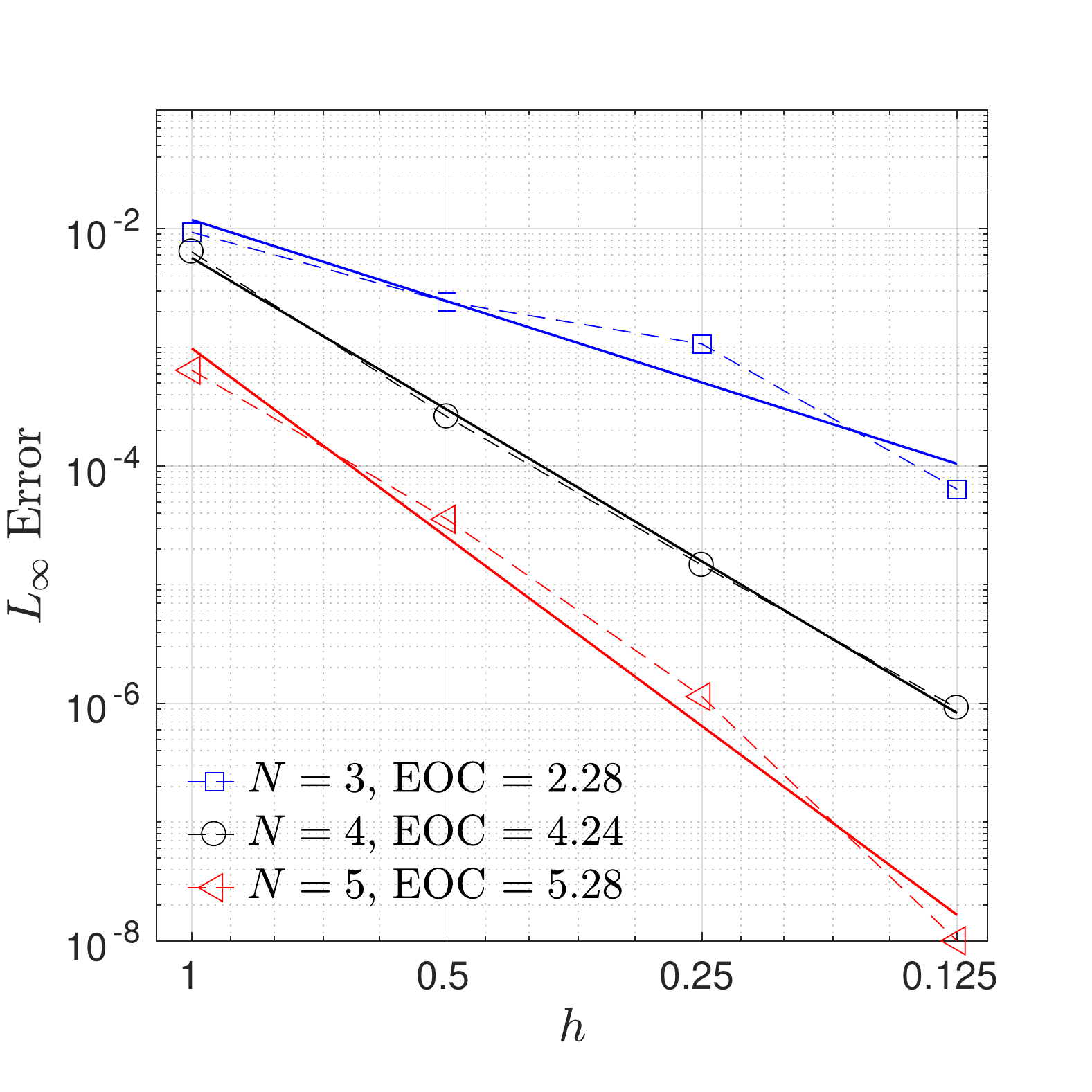}
		       \end{subfigure}
		       ~
		       \begin{subfigure}[b]{0.3\textwidth}
				\includegraphics[width=\textwidth]{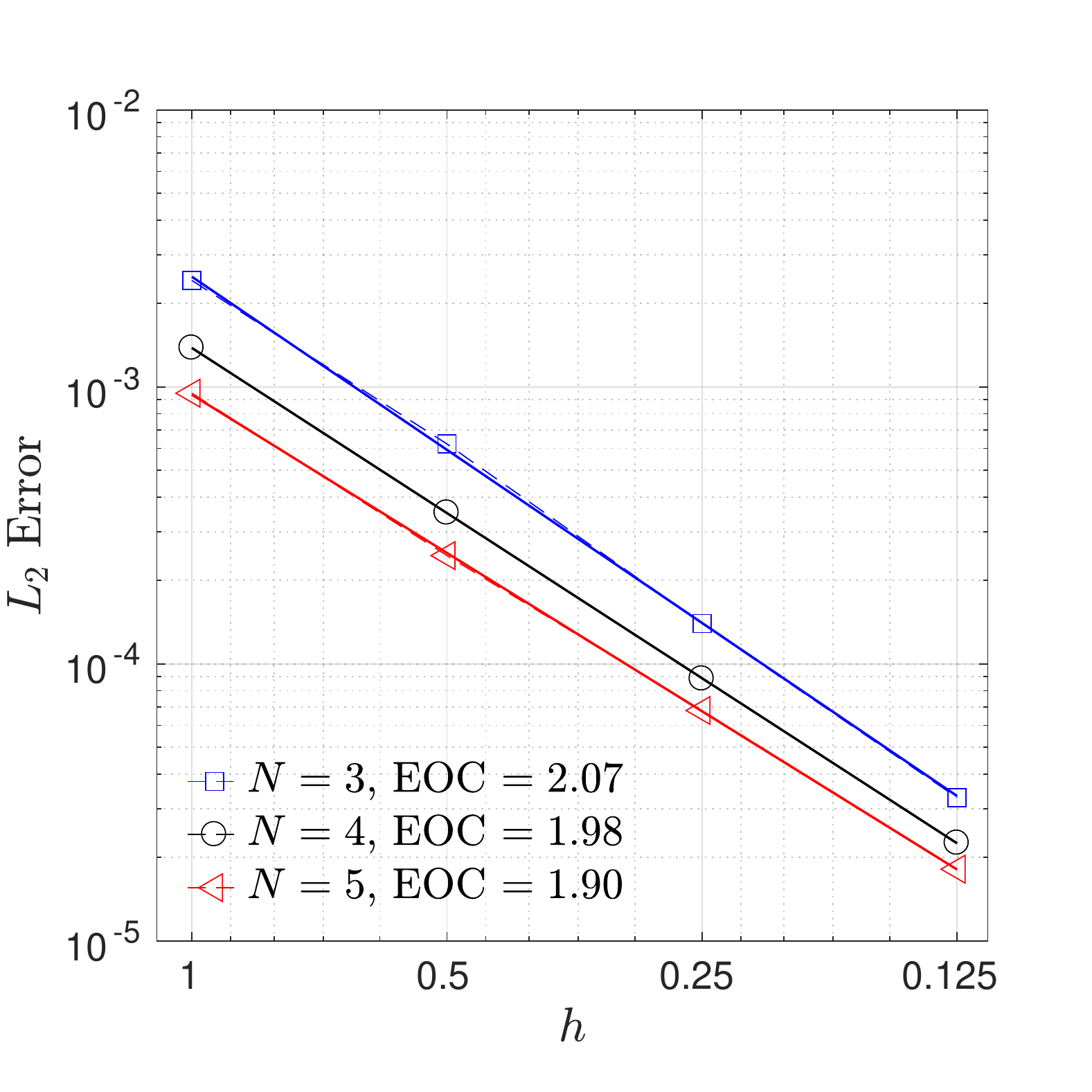}
		\end{subfigure}
		        ~
		        \begin{subfigure}[b]{0.3\textwidth}
		           \includegraphics[width=\textwidth]{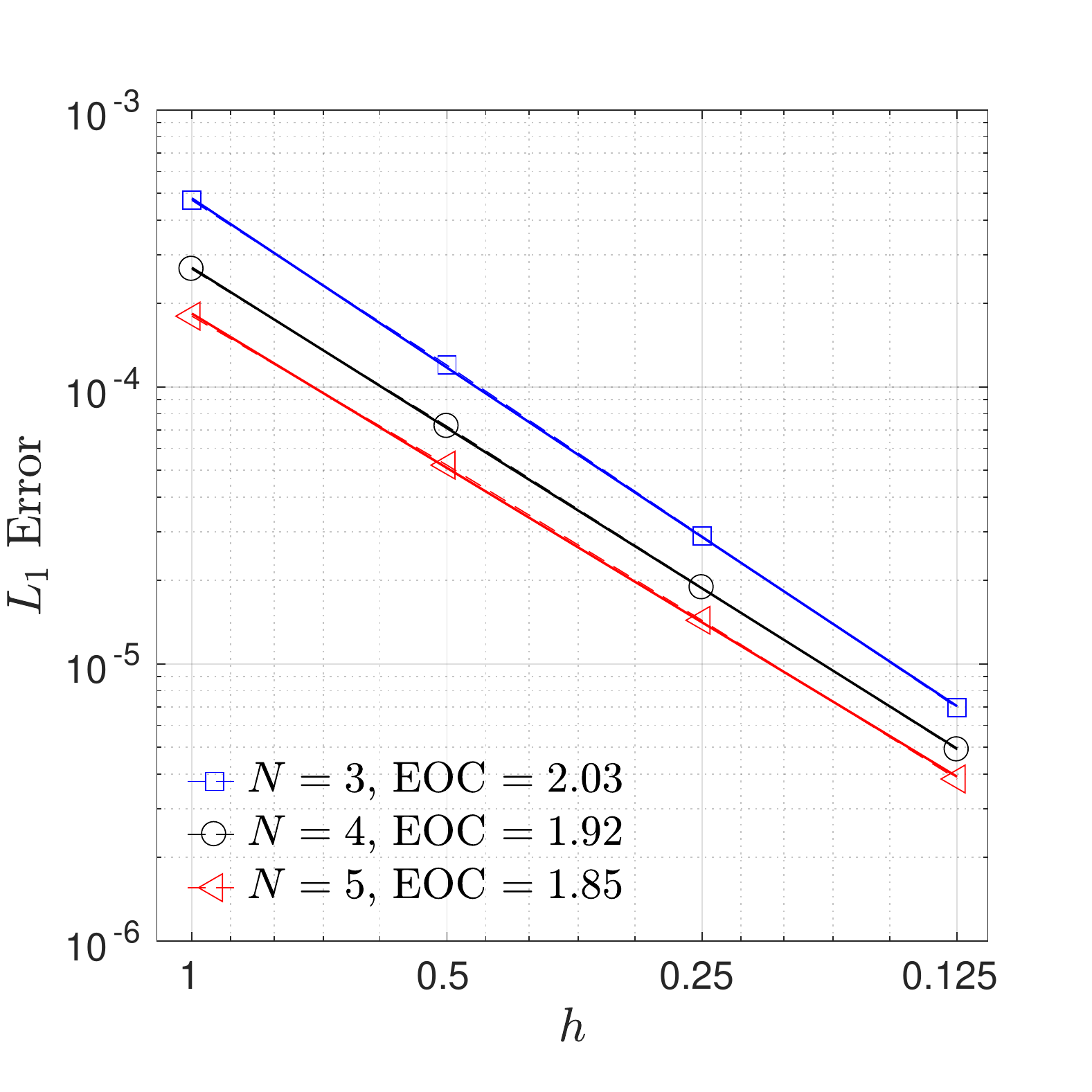}
		       \end{subfigure}
		        ~
		        \begin{subfigure}[b]{0.3\textwidth}
		           \includegraphics[width=\textwidth]{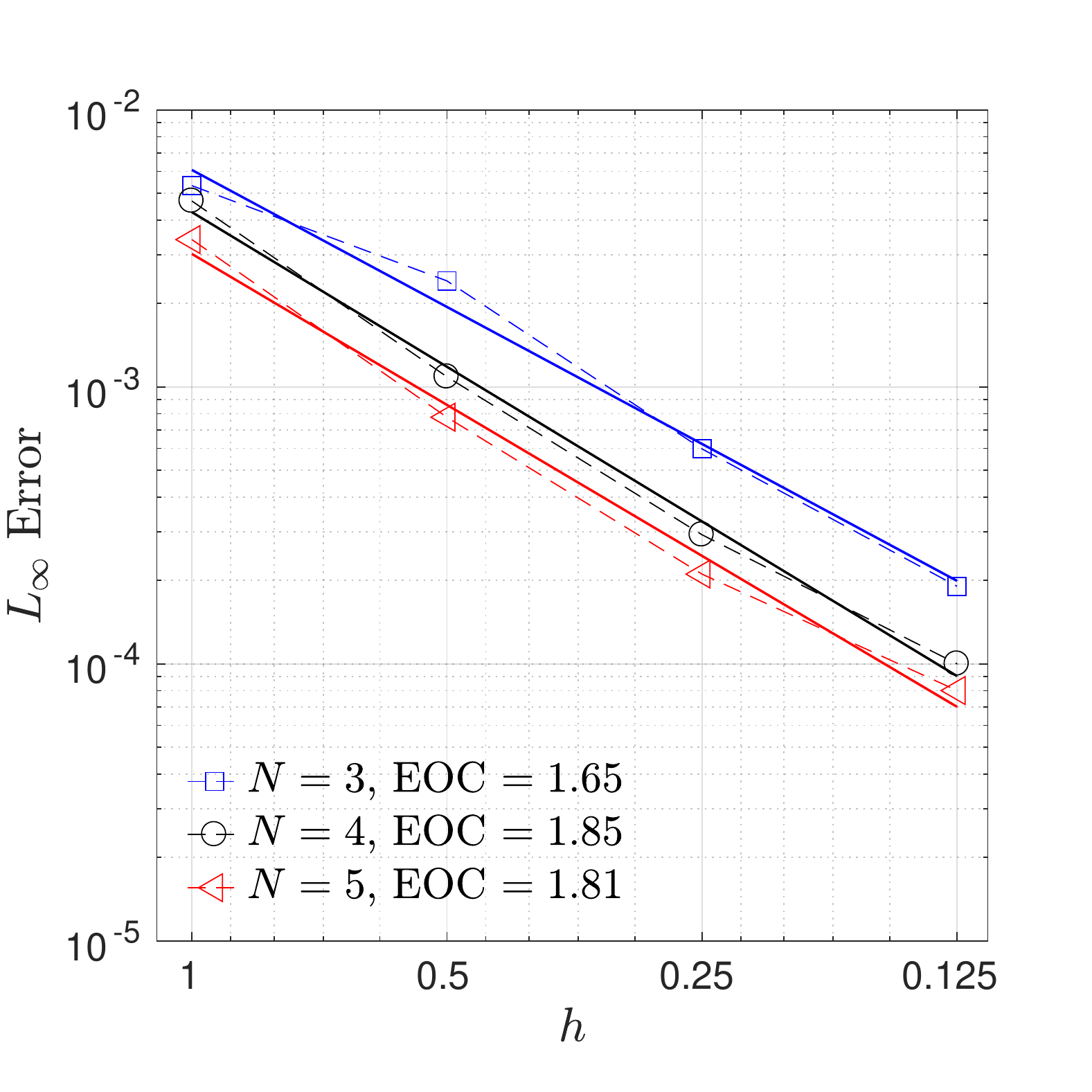}
		       \end{subfigure}
	\end{center}
	\caption{Spatial accuracy test for reinitialization of the level set function on the circular interface problem. $L_2$ and $L_{\infty}$ errors are computed on the band thickness of $\epsilon=0.3$ so that the kink point at the center of circle is avoided. $L_1$ error computation uses corresponding characteristic mesh length to smooth the Heaviside function. }
	\label{Fig.CircleTestConvergence}
\end{figure}

The first row of Figure \ref{Fig.CircleTestConvergence} shows the measured order of accuracy for $ N=3,4,5 $ using $L_2$, $L_{1}$ and $L_{\infty} $ error estimates near the interface with the condition of $ \epsilon =0.3 $ such that kink point at $ (0,0) $ is excluded. Although limiting is active for all numerical tests, troubled cell detector either do not mark elements in the narrow band region, or only a small number of elements are detected when the the solution is smooth and the mesh is fine enough to resolve high gradients near the interface.  Because of this, the numerical scheme reaches optimal $ N+1$  convergence rates near the smooth interface region. 

To illustrate the effect of subcell limiting on the convergence rates, which is not straight forward as the number of troubled elements are problem dependent, we artificially activated the subcell limiting on all cells for the same set of numerical experiments. The resulting convergence properties are given in the second row of Figure \ref{Fig.CircleTestConvergence}. Our scheme achieves its designed second-order accuracy  in all error norms.

\begin{figure}[ht!]
	\begin{center}
		\begin{subfigure}[b]{0.32\textwidth}
				\includegraphics[width=\textwidth]{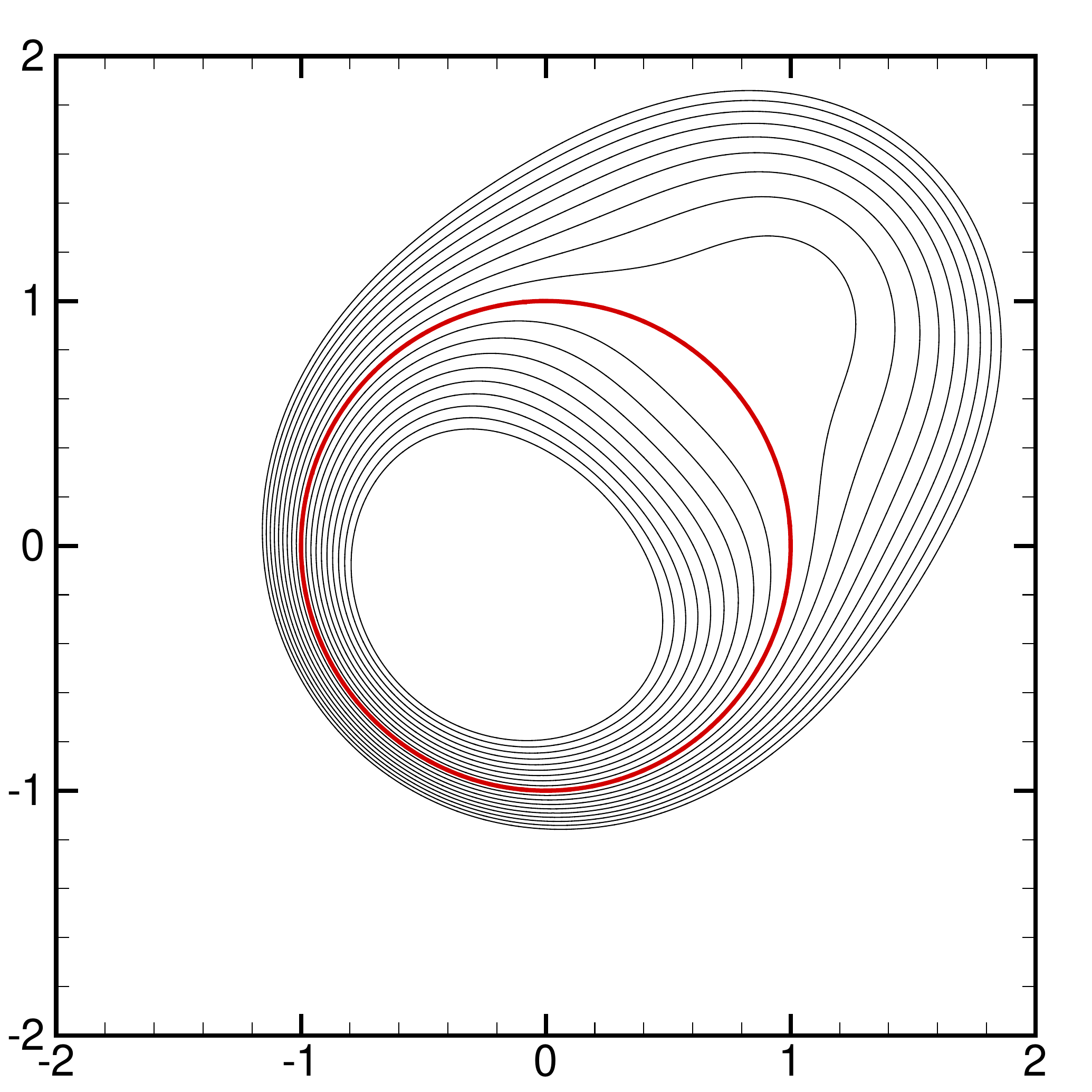}    
				\caption{$t=0s $}  
		\end{subfigure}
		  ~
		        \begin{subfigure}[b]{0.32\textwidth}
		           \includegraphics[width=\textwidth]{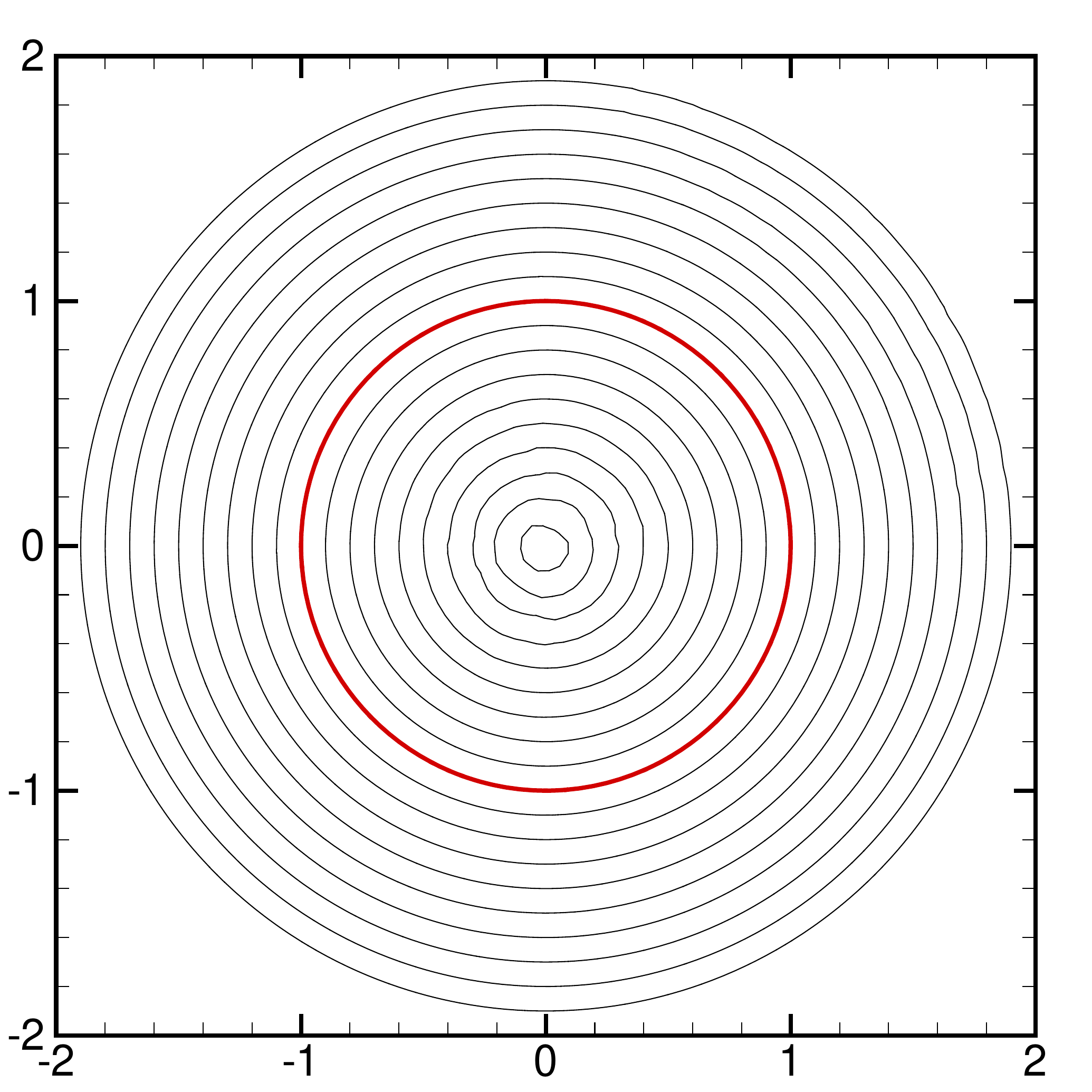}
		          \caption{$h $} 
		       \end{subfigure}
		        ~
		        \begin{subfigure}[b]{0.32\textwidth}
		           \includegraphics[width=\textwidth]{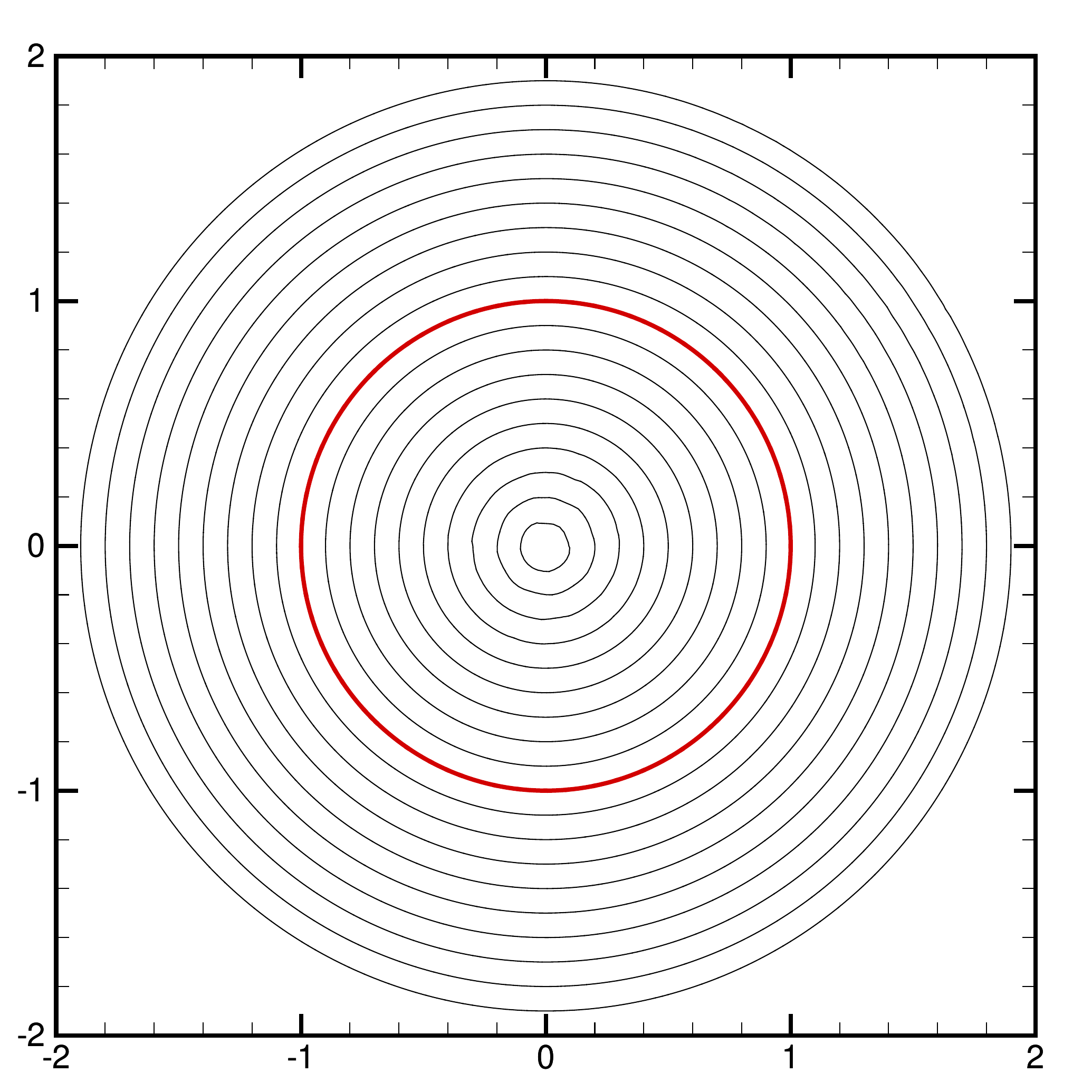}
		          \caption{$h/2 $} 
		       \end{subfigure}
	\end{center}
	\caption{Reinitialization of the level set function for the circle test on grids with characteristic length $ h $, $h/2$ and $h/4$ for $ N=5 $. Drawn are contour levels from $ -0.9 $ to $ 0.9 $ with step size  $ 0.1 $.}
	\label{Fig.CircleTestEvolution}
\end{figure}

Figure \ref{Fig.CircleTestEvolution} illustrates the reinitialized level set function for the computational grids of $h$, $h/2$, and $h/4$ for the approximation order $N=5$. A good recovery of the signed distance function is obtained for highly disturbed initial data even on the coarse grid. Upon increasing the resolution, the accuracy of the solution is improved around the kink point as expected. We also emphasize that proposed method does not cause the artificial movement of the interface after reinitialization as supported by the computed $L_1$ errors in the first row. 

As a second smooth interface problem, we consider the computation of signed distance function  from the reinitialization of an ellipse starting with the following initial LS function,
\begin{equation*}
\phi_{0}\left(x,y\right)=\left(\left(x-x_{0}\right)^{2}+\left(y-y_{0}\right)^{2}+0.1\right) \left(\sqrt{\frac{x^{2}}{A^2}+\frac{y^{2}}{B^2}}-1\right)
\end{equation*}
with $ A=1 $,$ B=0.5 $, $ x_{0}=0.875 $ and $ y_{0}=0.5 $. Similar to the circle test, the initial LS has both small and large gradients near the interface with an extended kink region between the line segment $ (-A,A) $ on the $x$ axis. In addition, the interface has a non-smooth region at the intersection point with  $y$ axis.  The closed form of the exact signed distance function is not known for the ellipse so it is accurately approximated by creating a finite number of points on the interface with coordinates,
\begin{equation*}
x_{n}=A\cos\left(2\pi n/N_n\right)\quad \text{and} \quad y_{p}=B\sin\left(2\pi n/N_n\right)
\end{equation*} 
where $ n $ and $ N_{n} $ are the index and total number of the points  inserted on the interface respectively. Then, the approximated  distance function is defined as
\begin{equation*}
\phi^{exact}(x_{i},y_{i})=\min(\sqrt{ (x_{i}-x_{n})^{2} +(y_{i}-y_{n})^{2}} ) \text{sgn}(\phi_{0} (x_{i},y_{i})).
\end{equation*}

\begin{figure}[ht!]
	\begin{center}
		\begin{subfigure}[b]{0.32\textwidth}
				\includegraphics[width=\textwidth]{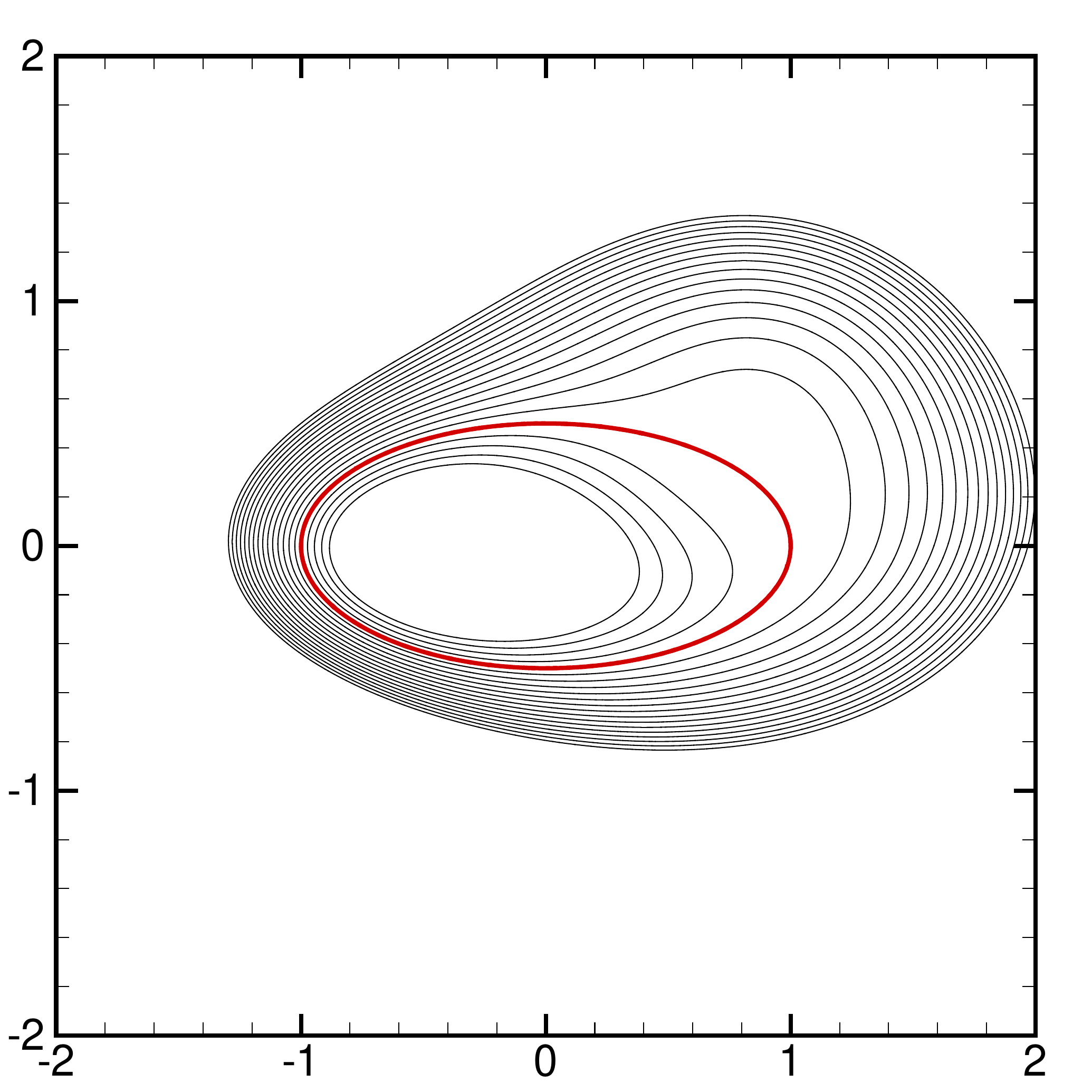}    
				\caption{$t=0s $}  
		\end{subfigure}
		  ~
		        \begin{subfigure}[b]{0.32\textwidth}
		           \includegraphics[width=\textwidth]{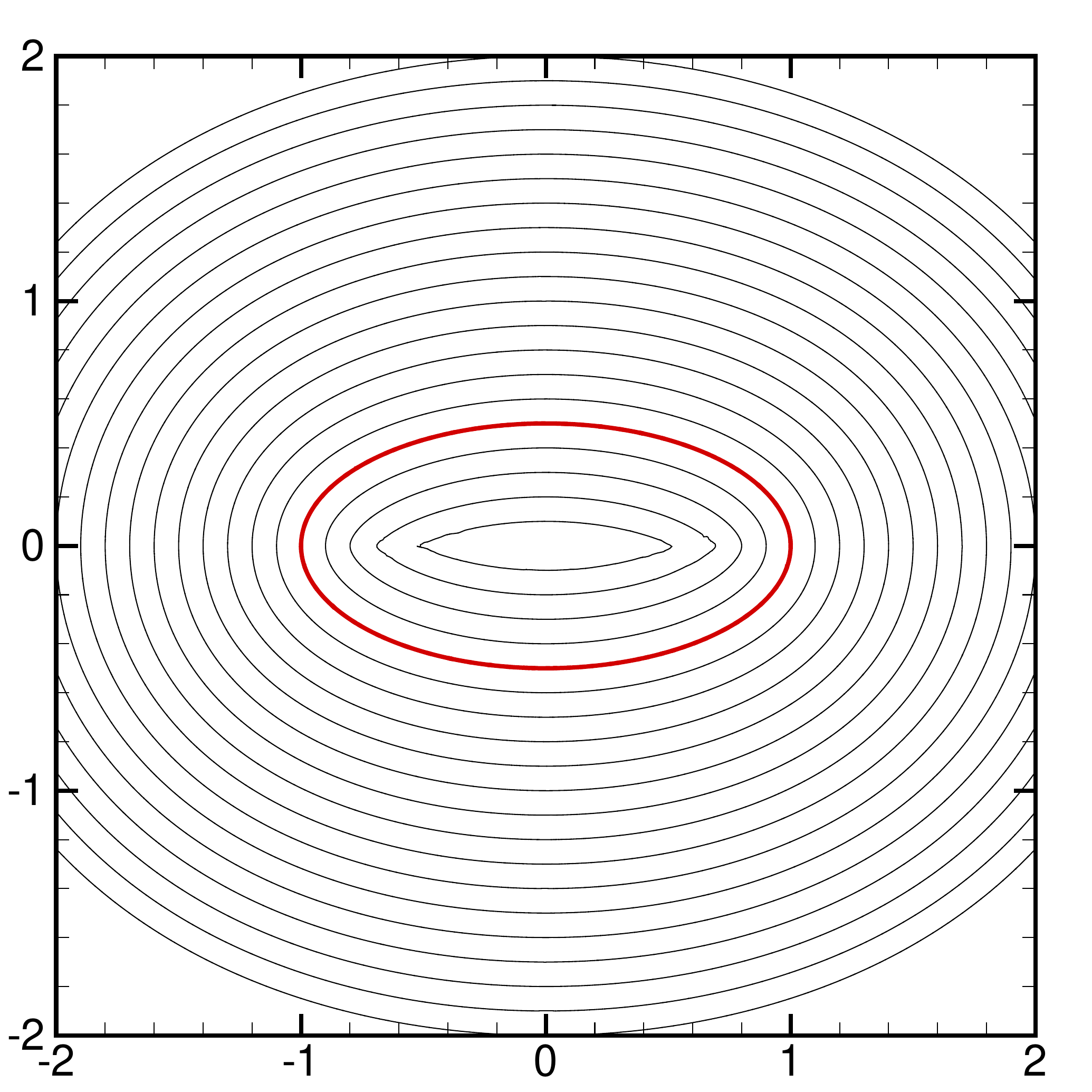}
		          \caption{$h/2 $} 
		       \end{subfigure}
		        ~
		        \begin{subfigure}[b]{0.32\textwidth}
		           \includegraphics[width=\textwidth]{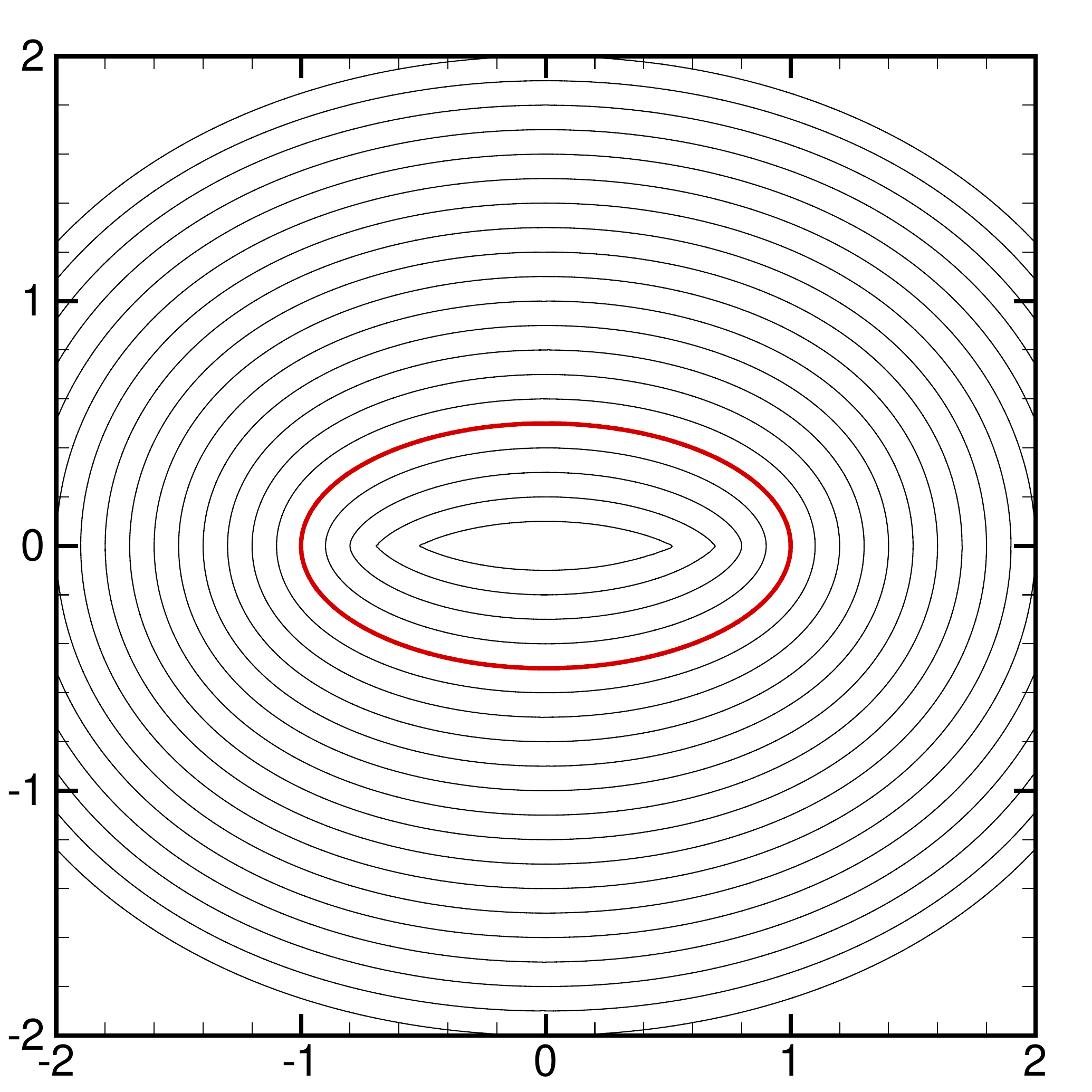}
		          \caption{$h/4 $} 
		       \end{subfigure}
	\end{center}
	\caption{Reinitialization of the level set function for the elliptic interface test for grids $ h/2 $, $h/4$ and for $ N=5 $ at initial state (a) and final time of $t=1.5$ (b-c). Drawn are contour levels from $ -0.4 $ to $ 1.5 $ with step size  $ 0.1 $.}
   \label{Fig.EllipseTestEvolution}
\end{figure}

In all tests, we use the same mesh configuration as with the circle test case above. Figure \ref{Fig.EllipseTestEvolution} illustrates the reinitialization of the highly perturbed LS function at different solution times for the meshes $ h/2 $ and  $ N=5$. As seen in the figure, the signed distance function is recovered well from the strongly distorted initial data.  

\begin{figure}[ht!]
	\begin{center}
		\begin{subfigure}[b]{0.4\textwidth}
				\includegraphics[width=\textwidth]{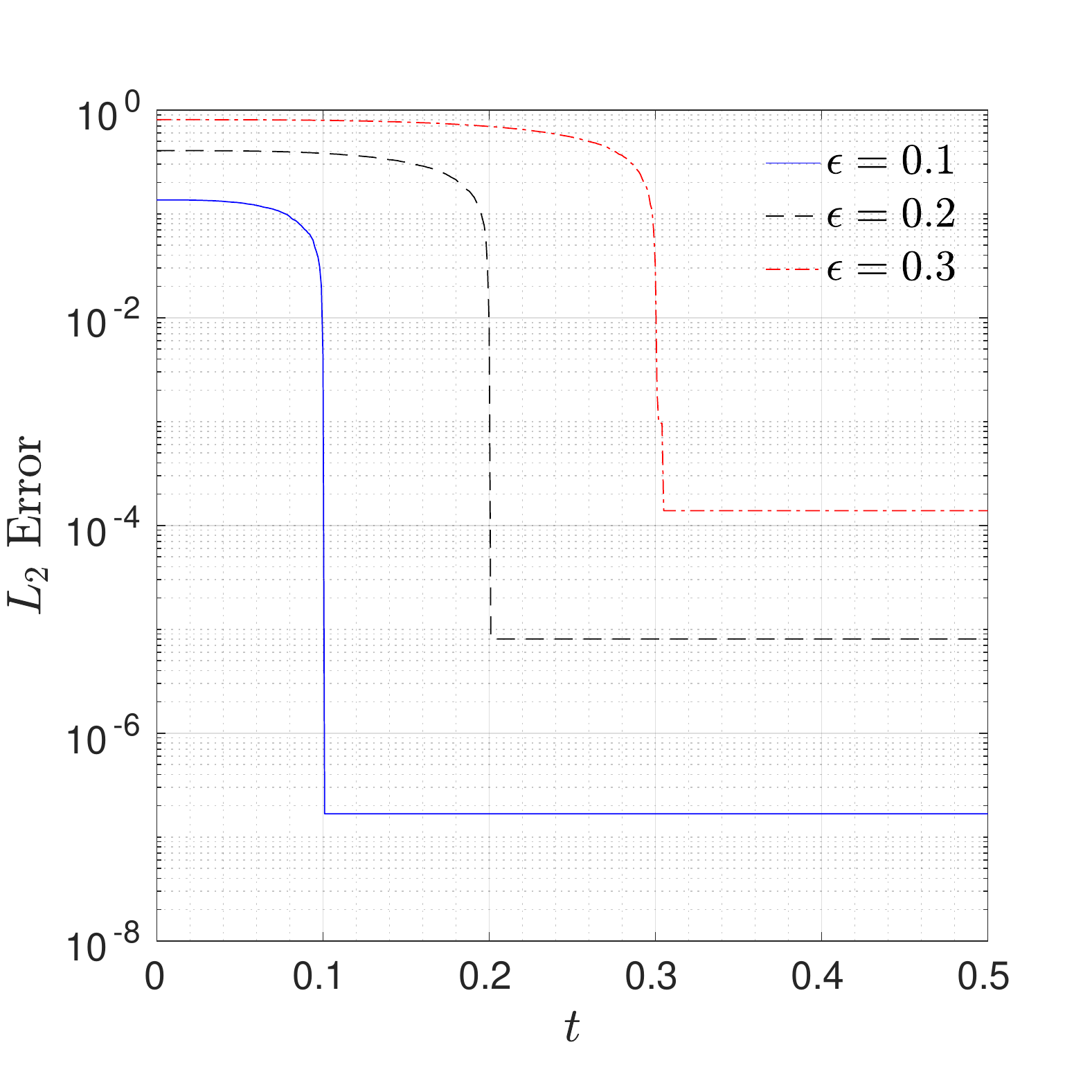}    
				\caption{}  
		\end{subfigure}
		  ~
		        \begin{subfigure}[b]{0.4\textwidth}
		           \includegraphics[width=\textwidth]{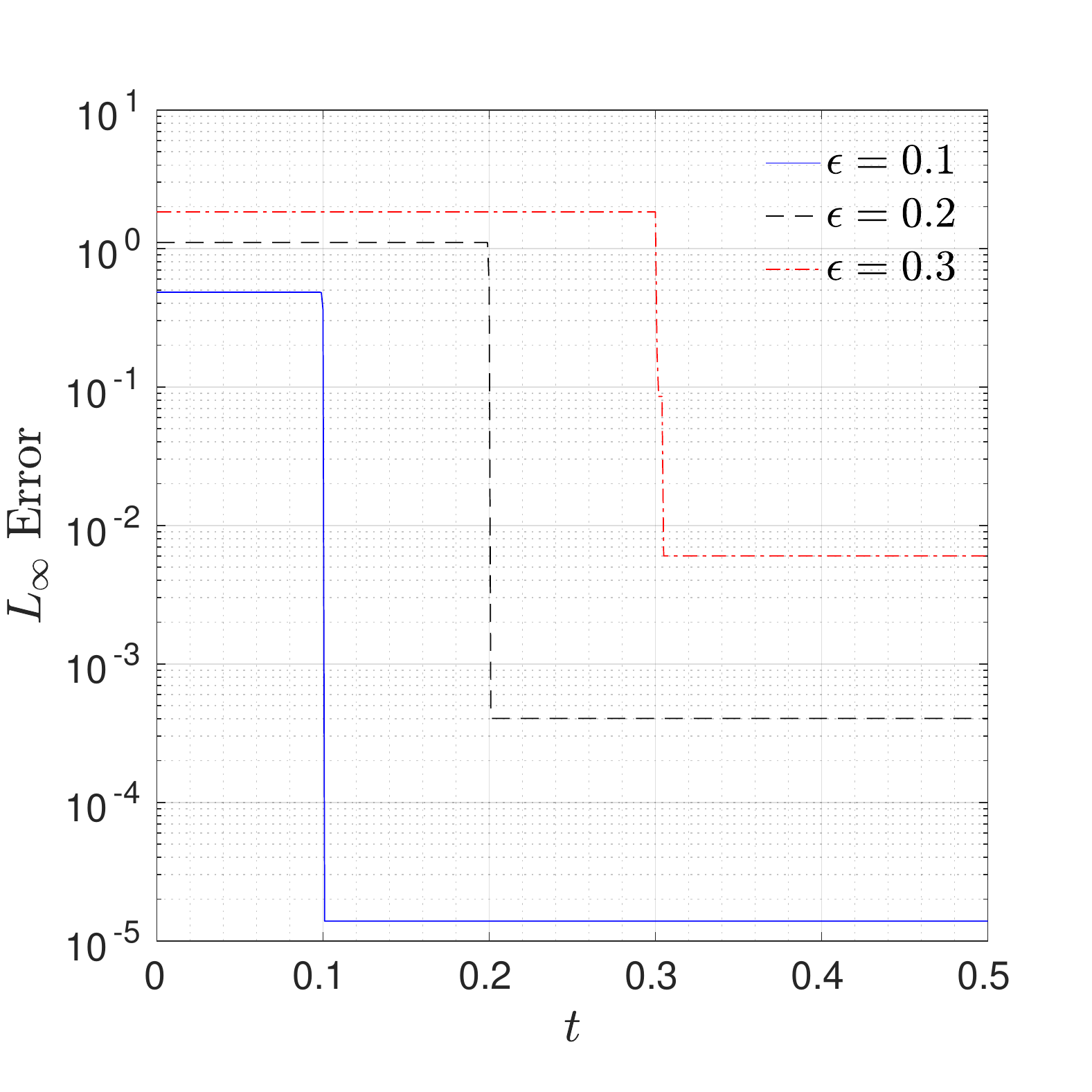}
		          \caption{} 
		       \end{subfigure}
		      ~
	\end{center}
	\caption{Local level set reinitialization for the elliptic interface test for grid $ h/2 $ and $ N=5 $. Time history of $L_2$ (a) and $L_\infty$ (b) norms of error measured on band thicknesses of $\epsilon = 0.1,0.2,03$. }
	\label{Fig.EllipseLocal}
\end{figure}
To illustrate the effectiveness and accuracy of the reinitialization method in local problems, we measure the error computed on a narrow band using three different thicknesses,$\epsilon=0.1, 0.2,0.3$. Figure \ref{Fig.EllipseLocal} shows the change of $L_2$ and $L_\infty$ norms of error in time on the relatively coarse grid $h/2$ and polynomial order $N=5$. For each local level set representation, the signed distance function is obtained at around the same time with the band thicknesses. With the increasing band thickness, error in both norms increases with including more troubled elements in error computations. Also, the test reveals that that pollution caused by the limiting does not spread out and degrade the solution outside of the region having kinks.

\subsection{Non-smooth Interfaces}
In the first non-smooth interface test, reinitialization of two intersecting circles of radii $ r $ centered at $ (\pm a,0) $ and $ 0<a<r $ is considered. Because $ 0<a<r $ holds, the circles intersect and the interface of interest is the union of the two circles. The signed distance function to the interface is given as
\begin{equation*}
\label{Eq.Test2}
\resizebox{0.99\hsize}{!}{$
\begin{aligned}
d(x,y)=
\begin{cases}
\min\left(\sqrt{x^2+\left(y \pm \sqrt{r^2-a^2}\right)^2}\right) &\ \; \text{if} \; \frac{a-x}{\sqrt{\left(a-x\right)^2 +y^2}} \geqslant \dfrac{a}{r} \enskip \text{and} \enskip \frac{a+x}{\sqrt{\left(a+x\right)^2 +y^2}} \geqslant \frac{a}{r} \\
\min\left(\sqrt{\left(x\pm a\right)^2+y^2}-r\right)    &\  \text{else}.\\
\end{cases}  
\end{aligned}
$}
\end{equation*}
Comparing with the previous two tests, this one is more critical because the signed distance function has kinks on  the whole $ y $ axis for $\phi>0$ and line segment $ [-a,a] $ on the $ x $ axis. The problem is solved for $ r=1 $ and $ a=0.7 $ on a computational domain of $ [-2,2]^2 $. Similar to the previous tests, the initial LS function is defined by multiplying the  signed distance with a perturbation function to create highly varying gradients near the interface as follows
\begin{equation*}
\label{Eq.Ch2.IntersectingCircles}
\phi_{0}\left(x,y\right)=\left(\left(x-1\right)^{2}+\left(y-1\right)^{2}+0.1\right)d\left(x,y\right).
\end{equation*}

\begin{figure}[ht!]
	\begin{center}
		\begin{subfigure}[b]{0.32\textwidth}
				\includegraphics[width=\textwidth]{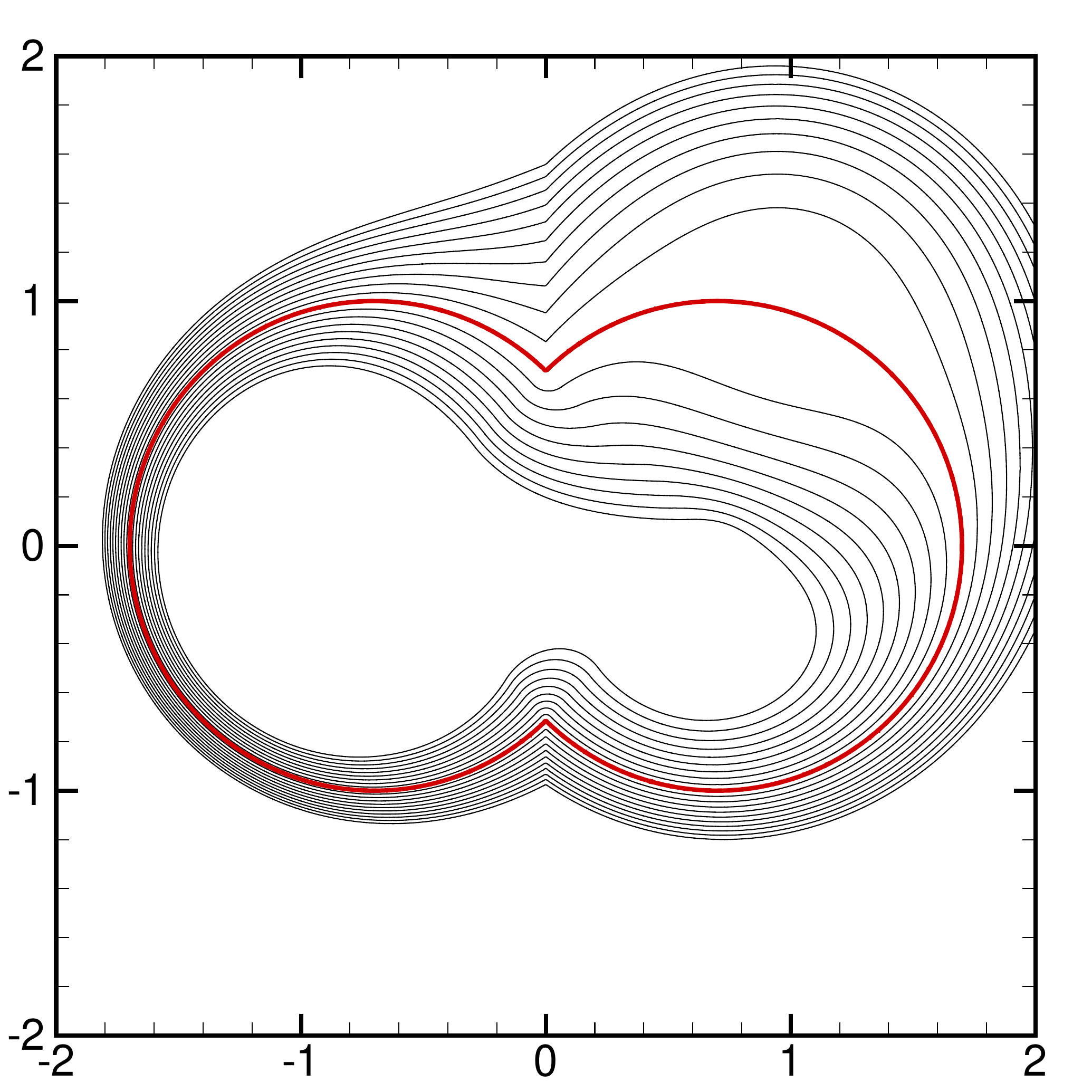}    
				\caption{$t=0s $}  
		\end{subfigure}
		  ~
		        \begin{subfigure}[b]{0.32\textwidth}
		           \includegraphics[width=\textwidth]{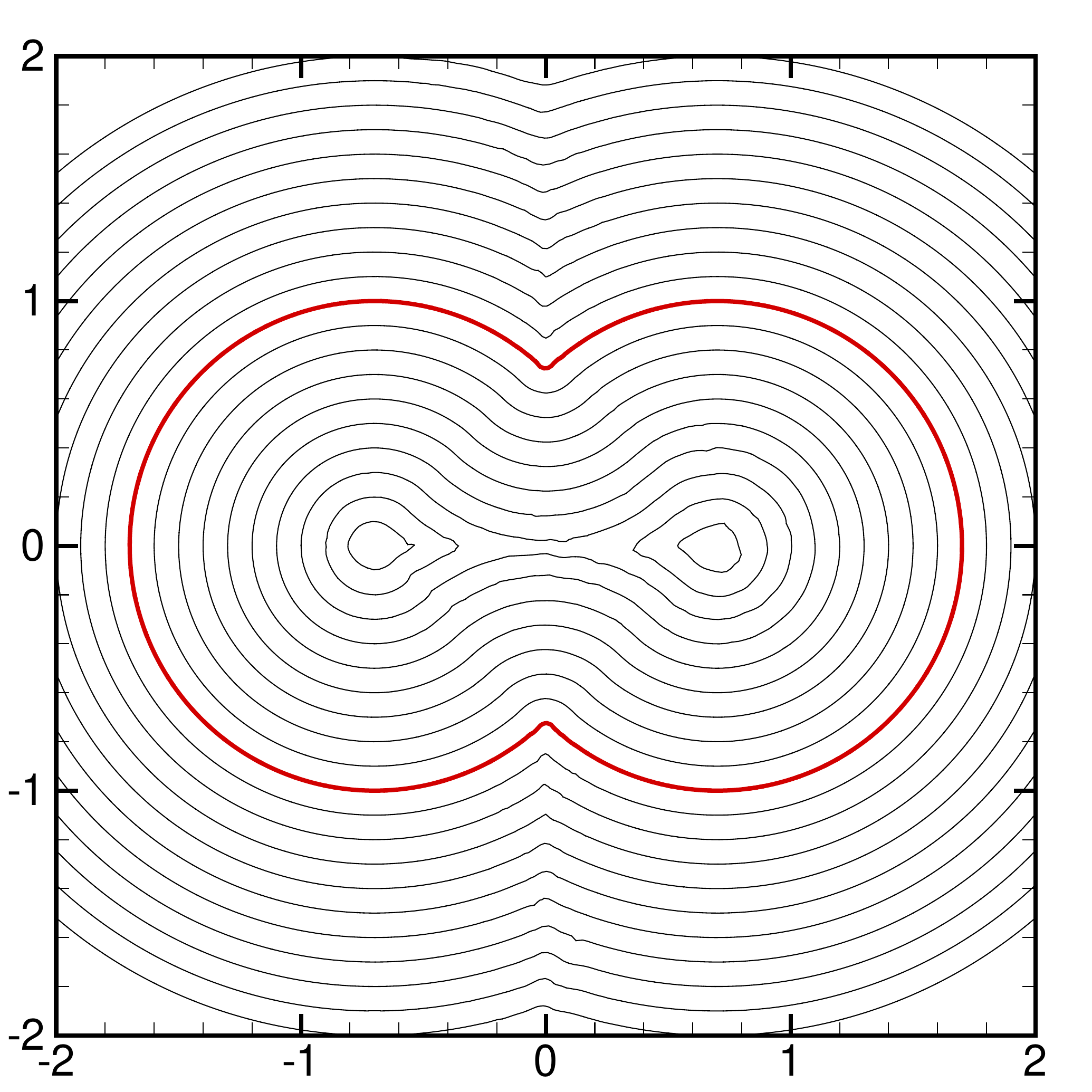}
		          \caption{$h/2 $} 
		       \end{subfigure}
		        ~
		        \begin{subfigure}[b]{0.32\textwidth}
		           \includegraphics[width=\textwidth]{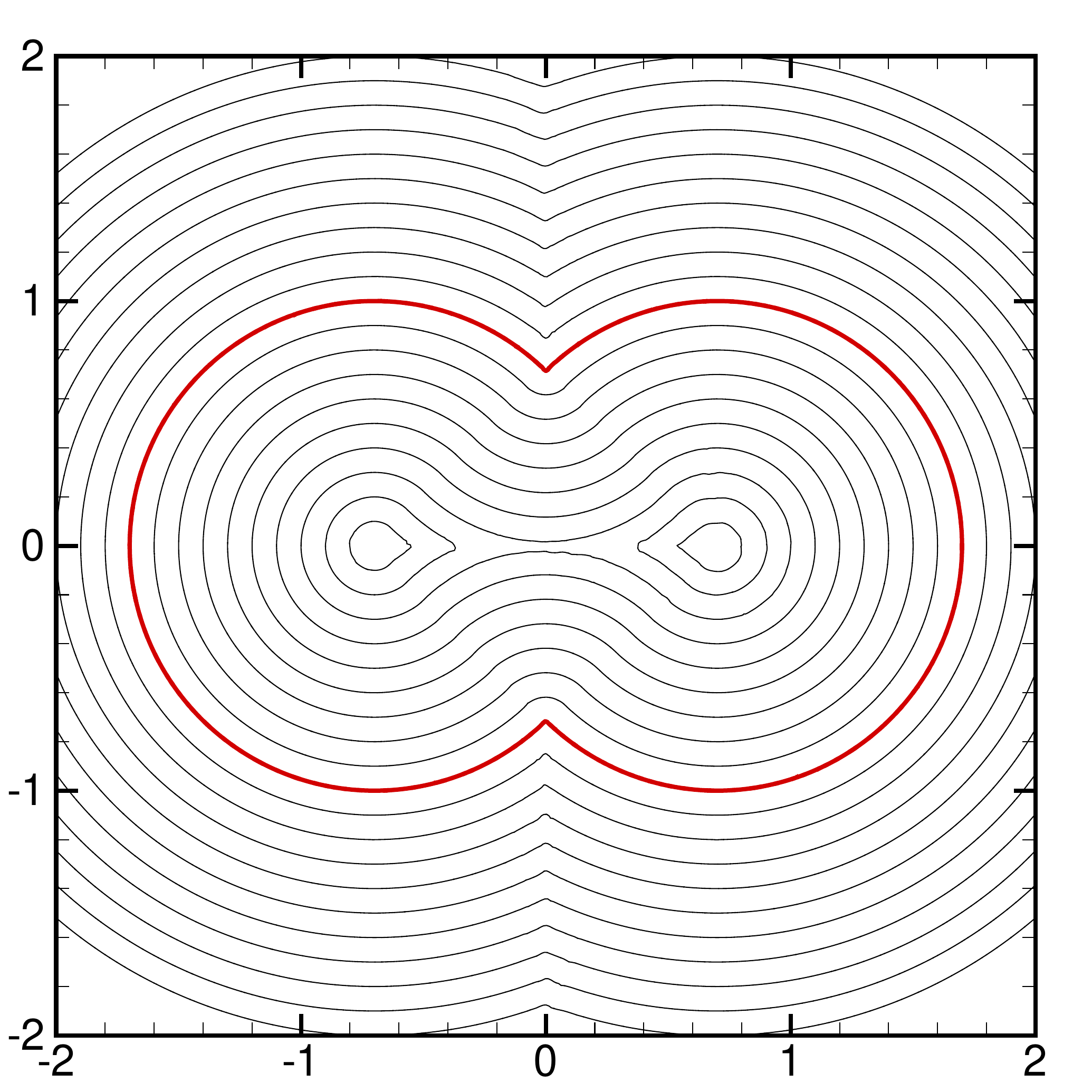}
		          \caption{$h/4 $} 
		       \end{subfigure}
	\end{center}
	\caption{Reinitialization of the level set function for the intersecting circles test for grids $ h/2 $, $h/4$ and  $ N=5 $ at initial state $t=0$ (a) and final time $t=1.0$ (b-c). Drawn are contour levels from $ -0.9 $ to $ 0.9 $ with step size  $ 0.1 $.}
	\label{Fig.IntersectingCircleTest}
\end{figure}

Figure \ref{Fig.IntersectingCircleTest} shows the reinitialization of  LS function at different solution times including the initial solution for the meshes $ h/2 $ and $h/4$, and  polynomial order $ N=5$. As seen in the figure, the signed distance function is generated accurately without any stability issue and  the kinks are more resolved with increasing resolution. 

As a second non-smooth test, construction of the signed distance function is considered for a square inside the computational domain of $ [-2,2]^2 $ starting with following initial LS function
\begin{equation*}
\phi_{0}(x,x)=0.8\max\left( \lvert x-x_{c} \lvert-w/2, \lvert y-y_{c} \lvert-w/2\right).
\end{equation*}

\begin{figure}[ht!]
	\begin{center}
		\begin{subfigure}[b]{0.32\textwidth}
				\includegraphics[width=\textwidth]{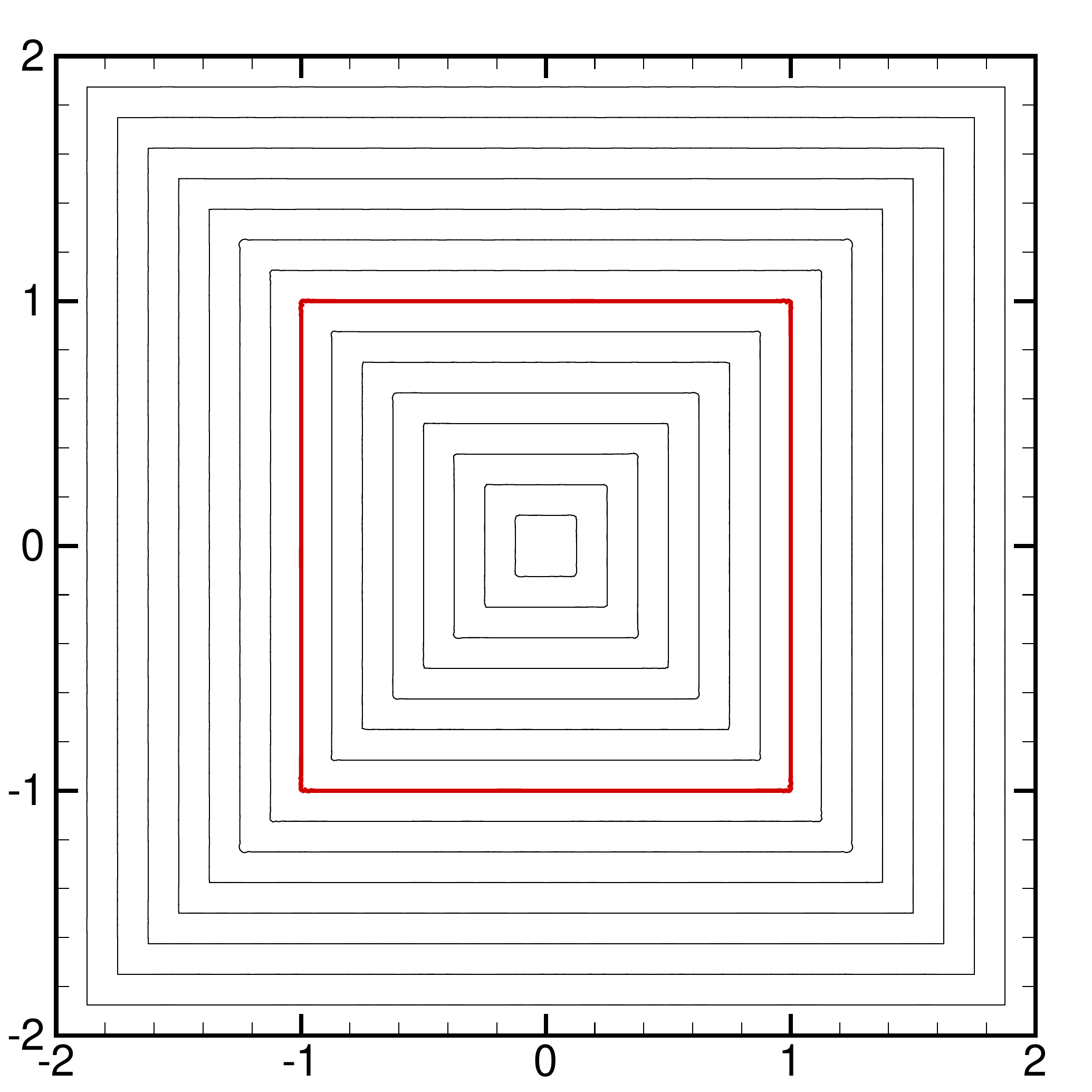}    
				\caption{$t=0s $}  
		\end{subfigure}
		  ~
		        \begin{subfigure}[b]{0.32\textwidth}
		           \includegraphics[width=\textwidth]{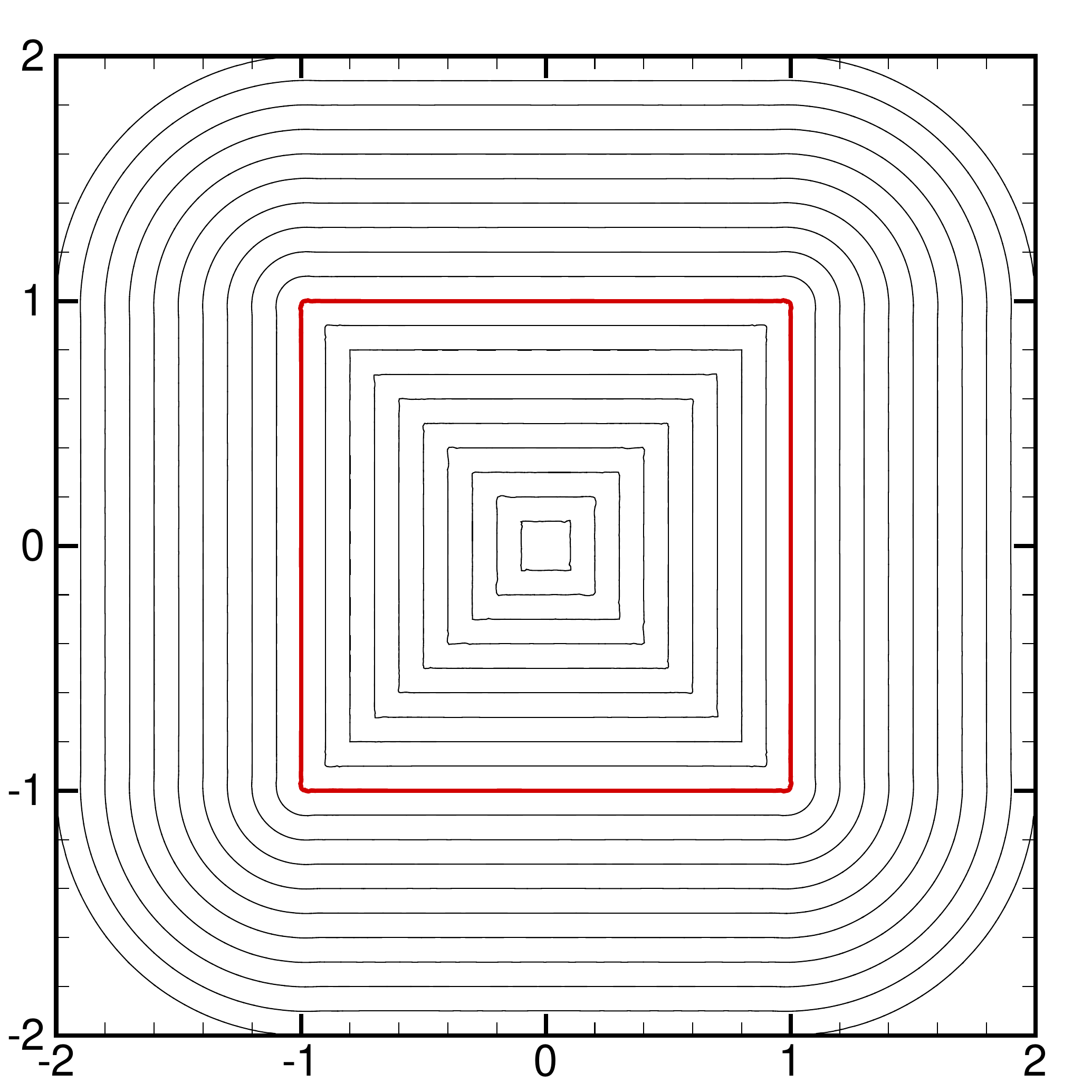}
		          \caption{$h/2 $} 
		       \end{subfigure}
		        ~
		        \begin{subfigure}[b]{0.32\textwidth}
		           \includegraphics[width=\textwidth]{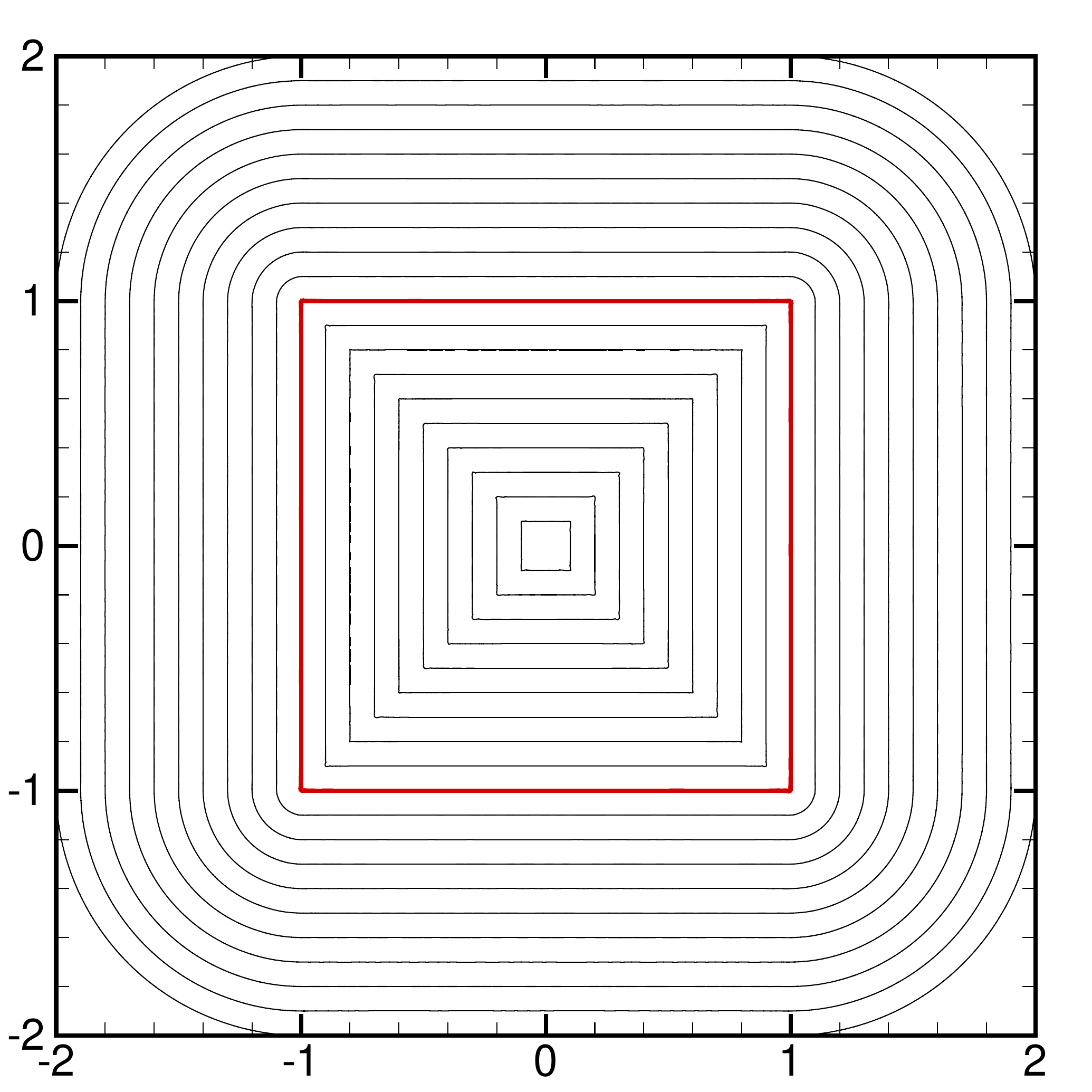}
		          \caption{$h/4 $} 
		       \end{subfigure}
	\end{center}
	\caption{Reinitialization of the level set function for the square interface test for grids $ h/2 $, $h/4$ and for $ N=5 $ at initial state (a) and final time (b-c),$t=1.5$. Drawn are contour levels from $ -0.9 $ to $ 1.0 $ with step size  $ 0.1 $.} \label{Fig.SquareTest}
\end{figure}

For $ w=2$, $ x_{c}=y_{c}=0$, the initial LS creates concentric squares centered at the origin with the interface having width of $ 2 $. Obviously, $\phi_{0} $ is not a signed distance function and includes kinks along the diagonals of the domain. The exact distance function is approximated similarly to the ellipse test and has kinks at the diagonals for $\phi \leq 0 $. The LS function is reinitialized well for smooth and non-smooth regions with sharp corners  as illustrated in the figure \ref{Fig.SquareTest}. 

\begin{figure}[ht!]
	\begin{center}
		\begin{subfigure}[b]{0.4\textwidth}
				\includegraphics[width=\textwidth]{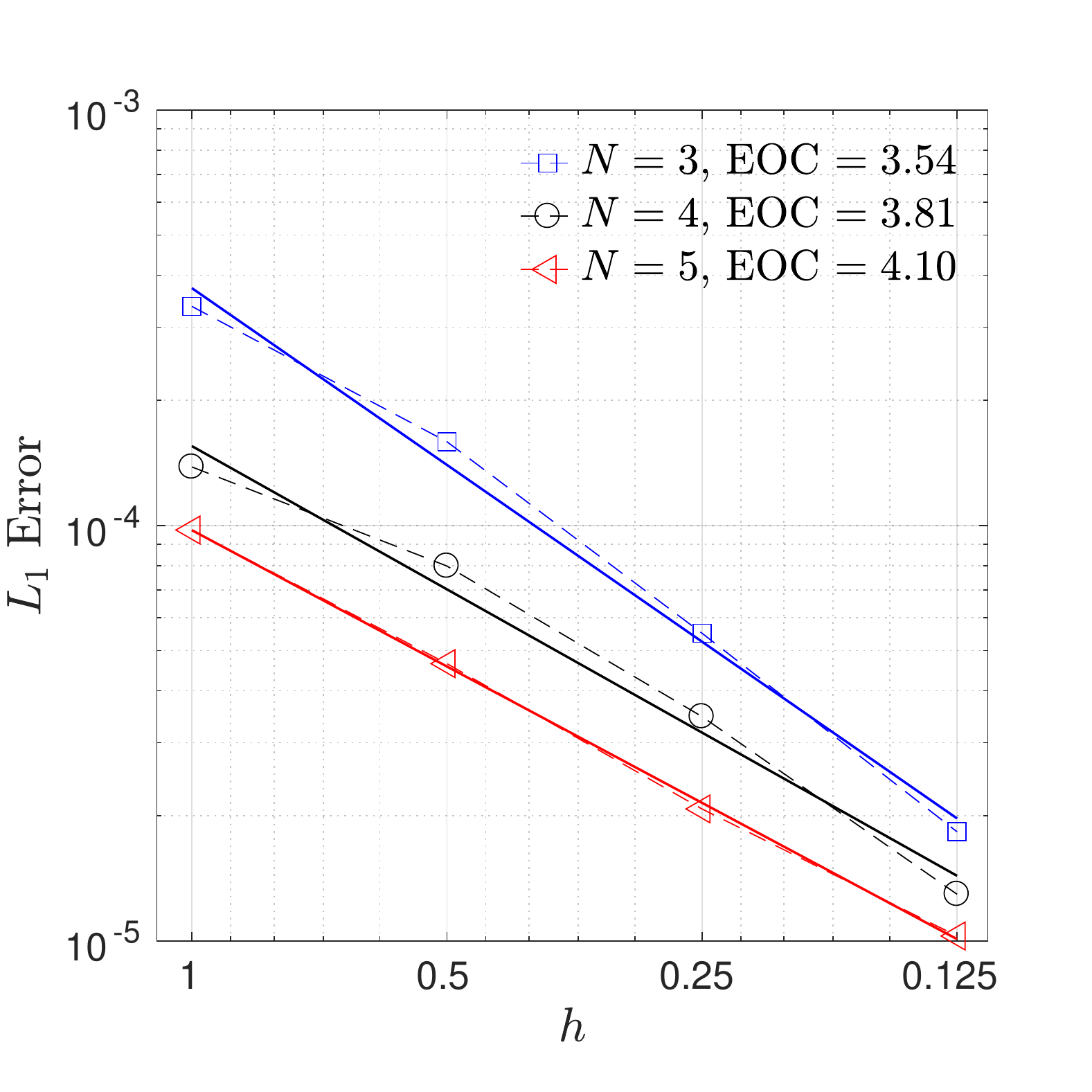}    
				\caption{Square Interface}  
		\end{subfigure}
		  ~
		        \begin{subfigure}[b]{0.4\textwidth}
		           \includegraphics[width=\textwidth]{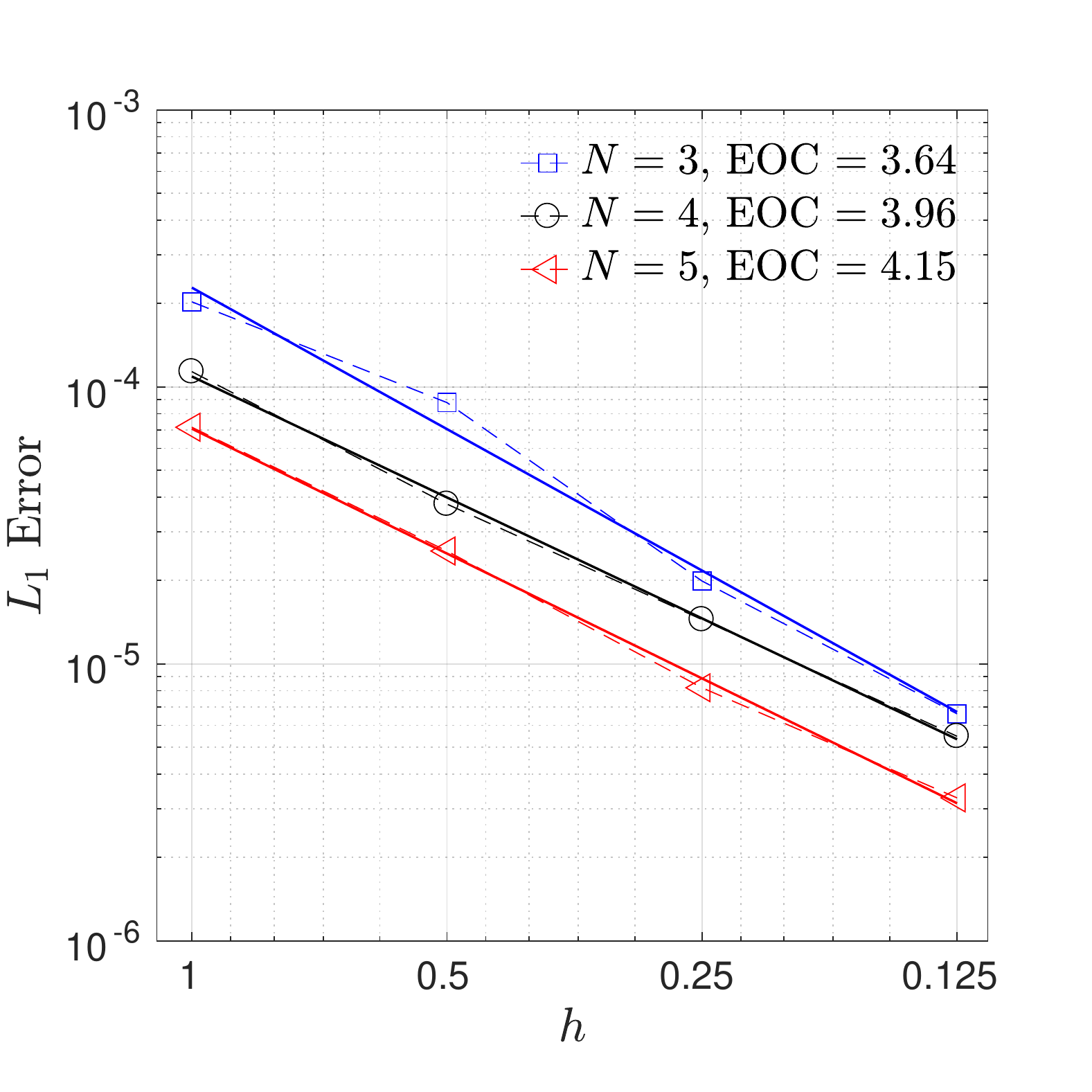}
		          \caption{Intersecting Circles} 
		       \end{subfigure}
	\end{center}
	\caption{Spatial accuracy test for non-smooth interface problems.   $L_1$ norm of errors are computed for square interface and intersecting circles tests.  }
	\label{Fig.L1NonSmooth}
\end{figure}
We emphasize that after projecting the initial condition to the approximation polynomial space, corners become rounded and lose their sharp profile. A better way of handling this deficiency is to refine the mesh around sharp corners and start with a more resolved initial condition \cite{karakus_gpu_2016}, but this is out of the scope of this study. Even without any special treatment of the interface singularities, the proposed technique provides an accurate interface representation as shown in Figure \ref{Fig.L1NonSmooth}. The $L_1$ norm of errors are computed for the intersecting circles and square interface test problems for $N=3,4,5$ on the same sequence of meshes in the Figure. Estimated orders of convergence on the error norm suggest convergence rates between $N+1$ and $2$ resulting from the stabilization on  elements located on the interface. As the polynomial order increases, oscillations of initial condition at corners and deficiency from the optimal convergence rate of $(N+1)$ increases.

\subsection{Multiple Interfaces}
A final numerical test contains a more complex interface structure that might be more relevant to practical physical problems such as motion of multiple bubbles in multiphase flows. For this case, we have defined $12$ circular interfaces of the same radius distributed over the computational domain $ [-2,2]^2 $.  The signed distance function to the interface is given as the minimum of distance functions of each circle
\begin{align*}
   d\left(x,y\right) = \min\left(d_i(x,y)\right), \quad d_i\left(x,y\right) = \sqrt{(x-x_i)^2 + (y-y_i)^2}-r_i, 
\end{align*}
where $\left(x_i, y_i\right)$ and $r_i$ are the center coordinates and radius of circle $i$ for $i=1\ldots 12$. To make the problem more challenging, initial value of the level set function is obtained by multiplying the distance function with the following function similar to the previous test cases.  
\begin{equation*}
\phi_{0}\left(x,y\right)=\left(\left(x-1\right)^{2}+\left(y-1\right)^{2}+0.1\right)d\left(x,y\right)
\end{equation*}

\begin{figure}[ht!]
	\begin{center}
		\begin{subfigure}[b]{0.4\textwidth}
				\includegraphics[width=\textwidth]{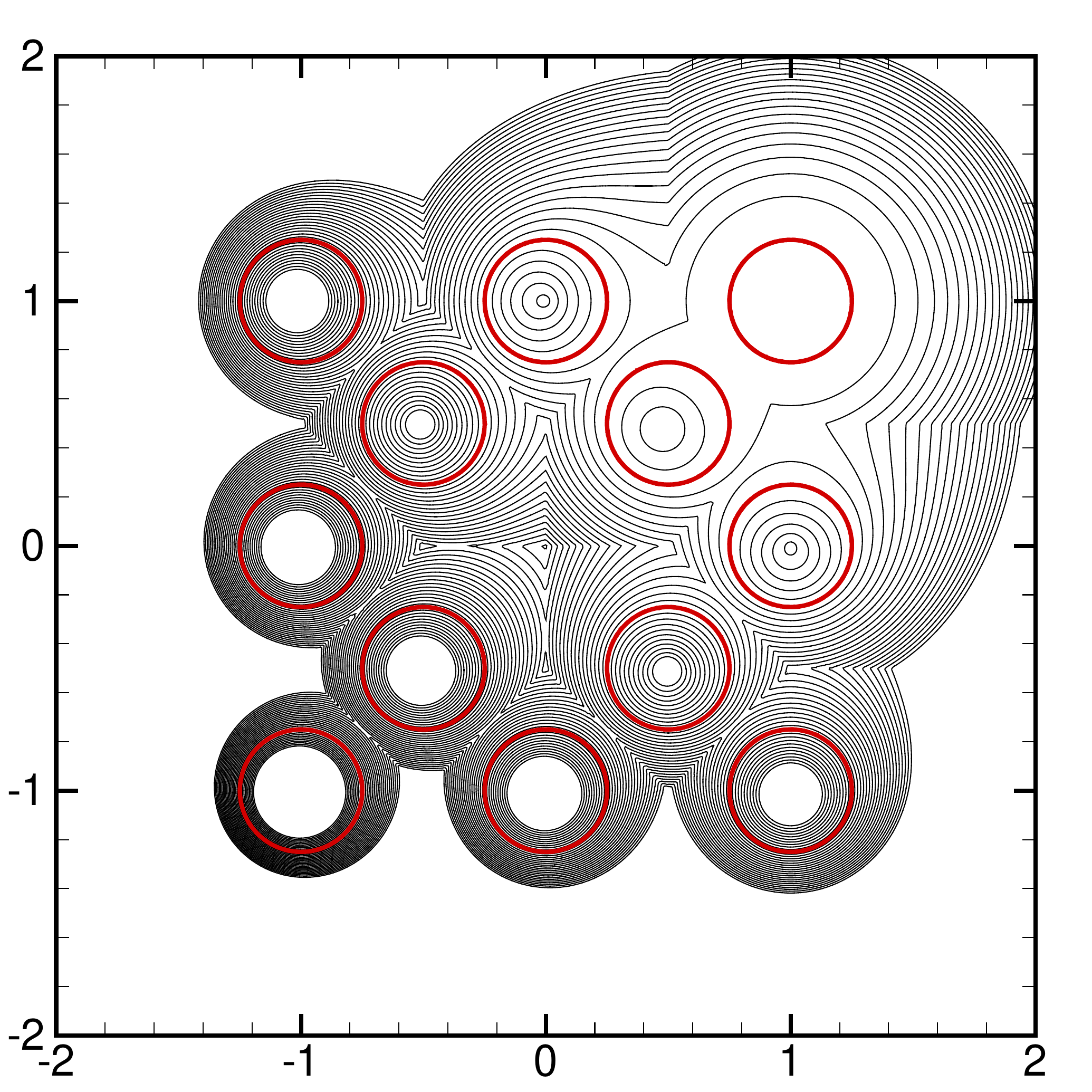}    
				\caption{$t=0s $}  
		\end{subfigure}
		  ~
		        \begin{subfigure}[b]{0.4\textwidth}
		           \includegraphics[width=\textwidth]{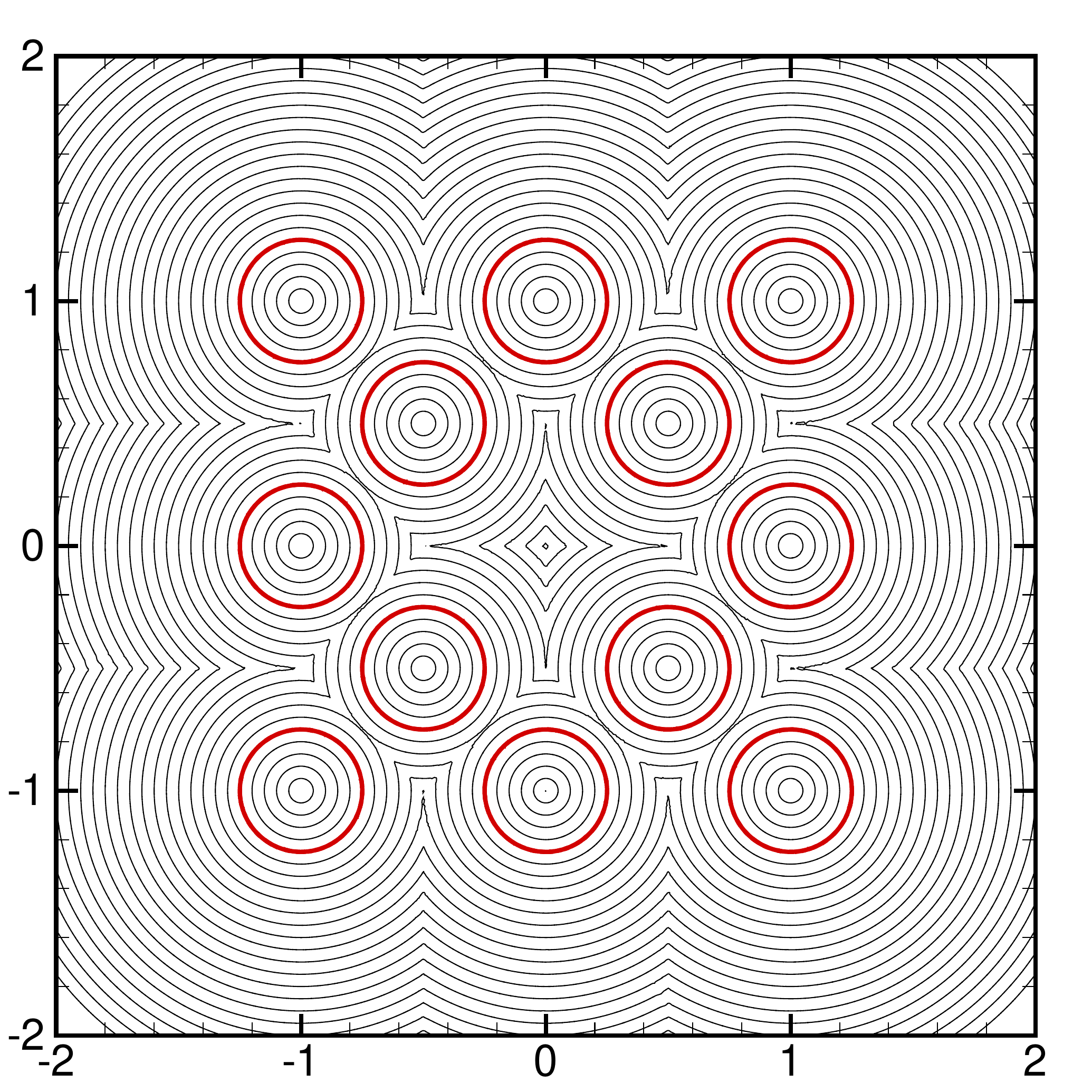}
		          \caption{grid $h/4 $} 
		       \end{subfigure}
	\end{center}
	\caption{Reinitialization of the level set function for the multiple interfaces test for the $ h/4 $ grid and $ N=5 $ at initial state (a), and final time, $t=1.1$ (b). Drawn are contour levels from $ -0.5 $ to $ 1.0 $ with step size  $ 0.05 $.}
	\label{Fig.MultipleInterfaces}
\end{figure}

The initial field has highly varying gradients and curvatures, and creates a complex structure having a range of kinks over the domain as shown in the Figure \ref{Fig.MultipleInterfaces} (a). We use this test case to illustrate stability of the algorithm in a more realistic case.  In Figure \ref{Fig.MultipleInterfaces}, the solution is obtained on the $h/4$ grid for $N=5$ with second-order finite volume subcell stabilization. The signed distance function is recovered well even under highly complex kink structures. 

%% file: conclusion.tex
We have presented a  high-order, local  discontinuous Galerkin approach to reinitialize level set functions through flow of the time Eikonal equation. To stabilize the resulting Hamilton-Jacobi equations, we utilized a subcell finite volume limiter based on second-order WENO reconstruction on triangular grids. We showed that the scheme achieves designed order of accuracy and preserves stability using smooth and non-smooth interface problems with highly varying gradients. The presented approach offers a high-order level set reinitialization by addressing stability and accuracy issues of standard  hyperbolic reinitialization in the DG framework \cite{adams_high-order_2019, utz_interface-preserving_2017}.

The presented solver is implemented on the open source project libParanumal (LIBrary of PARAllel NUMerical ALgorithms) \cite{ChalmersKarakusAustinSwirydowiczWarburton2020}. libParanumal consists of a collection of mini-apps with high-performance portable implementations of high-order finite element discretizations for a range of different fluid flow models \cite{karakus_discontinuous_2019,karakus_gpu_2019,swirydowicz_acceleration_2019}. 

The GPU performance of the scheme in triangular/tetrahedral elements remains to be investigated. Furthermore, extension to contact line problems and multiphase flows with sharp interfaces by considering potential performance gains  will be studied in future works.
